\documentclass{elsart}
\usepackage{stmaryrd}
\usepackage{amsfonts}
\usepackage{amsmath}
\usepackage{mathrsfs}
\usepackage{xcolor}
\usepackage[all]{xy}
\usepackage{graphicx}
\def\uuar{\mathord{\mbox{\makebox[0pt][l]{\raisebox{.4ex}
{$\uparrow$}}$\uparrow$}}}

\def\dda{\mathord{\mbox{\makebox[0pt][l]{\raisebox{-.4ex}
{$\downarrow$}}$\downarrow$}}}
\allowdisplaybreaks[4]
\newtheorem{tm}{Theorem}[section]
\newtheorem{pn}[tm]{Proposition}
\newtheorem{lm}[tm]{Lemma}

\theoremstyle{definition}
\newtheorem{dn}[tm]{Definition}
\newtheorem{rk}[tm]{Remark}
\newtheorem{ex}[tm]{Example}

\journal{Fuzzy Sets and Systems}
\usepackage{bm}

\let\sss=\scriptscriptstyle

\begin{document}
	
\begin{frontmatter}
\title{Topological representations for frame-valued domains via $L$-sobriety}
\author{Guojun Wu$^{1,2}$, Wei Yao$^{1,2}$, Qingguo Li$^{3}$}
\address{$^{1}$School of Mathematics and Statistics, Nanjing University of Information Science and Technology, Nanjing, 210044, China\\
$^{2}$Applied Mathematics Center of Jiangsu Province, Nanjing University of Information Science and Technology, Nanjing, 210044, China\\
$^{3}$School of Mathematics, Hunan University, Changsha, 410082, China
}

\address{}
\date{}

\begin{abstract}
With a frame $L$ as the truth value table, we study the topological representations for frame-valued domains. We introduce the notions of locally super-compact $L$-topological space and  strong locally super-compact $L$-topological space. These notions allow us to successfully represent  continuous  $L$-dcpo's and algebraic  $L$-dcpo's
 via $L$-sobriety. By
means of Scott $L$-topology and specialization $L$-order, we   establish a categorical isomorphism between  the category of the continuous (resp., algebraic) $L$-dcpo's with Scott continuous maps and that of the locally super-compact  (resp., strong locally super-compact) $L$-sober spaces with  continuous maps.   As an application of the resulting topological representations, for a continuous  $L$-ordered set $P$,  we   obtain a categorical isomorphism between  the category of directed completions of  $P$ with Scott continuous maps and that of the  $L$-sobrifications of the Scott space $(P, \sigma_{L}(P))$ with  continuous maps.

\end{abstract}

\end{frontmatter}

\section{Introduction}
 Domain Theory \cite{Gierz},  developed from continuous lattices introduced by Scott \cite{Scott-72} in the  1970s as a denotational model for functional languages.  The mutual transformations and interactions of  the mathematical structures of orders and topologies  are  fundamental features  of domain theory. Sobriety, an important concept in both domain theory and non-Hausdorff topology \cite{Jean-2013}, plays a critical role in connecting topological and ordered structures (see \cite{XU-TopAPpp}). The topological representation of various domains is a significant research topic in domain theory, with sober spaces serving as essential tools for these representations. It is well-known that every injective $T_0$ space is sober,  and  Scott established a categorical isomorphism between   injective $T_0$ spaces and  continuous lattices \cite{Scott-72}.  The topological representations of  continuous dcpo's via sobriety is discussed earlier in \cite{Banaschewski,Erne,R. Hofmann,Lawson}. Based on these works,  Zhao \cite{Zhao-Parial} systematically compiled the  topological representations of   continuous dcpo's and  summarized the following equivalent statements:
 \begin{itemize}
\item $X$ is homeomorphic to the Scott space of a continuous dcpo;
\item $X$ is  both a sober space and   a  C-space;

 \item $X$ is   a sober space that has  an injective hull;

  \item $X$ is   a sober space whose  open set lattice is completely distributive.

 \end{itemize}
In \cite{XU-mao-form}, Xu and Mao studied the topological representations of various domains and obtained that  continuous dcpo's can be represented by locally super-compact sober spaces. This representation of  continuous dcpo's is both simple and elegant,  as it is achieved through set-theoretic sequences of proof without  using  higher-order definitions or properties,  such as the complete distributivity of open set lattices and the categorical properties of topological spaces.

Quantitative domains \cite{Smyth} are the extensions of classical domains by generalizing ordered relations to  more general structures,
such as enriched categories, generalized metrics and lattice-valued  orders.
The lattice-valued  order approach to quantitative domain theory, known as lattice-valued domain theory,  has been primarily developed by Wagner \cite{Wagner}, Fan and Zhang \cite{FanTCS,FanZhang},  Zhang et al. \cite{Zhang4,Yu-Zhang}, Yao et al. \cite{PartI,PartII}, Li et al. \cite{LiQG1,Su-Li-Intel}, Zhao et al. \cite{Zhao2,Zhao3,Wang-Zhao}, Yue et al. \cite{Liu-Meng-Yue}, Guti\'{e}rrez Garc\'{i}a et al.\cite{Garcia}, etc.
In \cite{YaoTFS}, Yao demonstrated  that the category of continuous $L$-lattices is isomorphic to that of injective $T_0$ $L$-topological spaces.  This indicates that   continuous $L$-lattices can be represented by  injective objects in $L$-{\bf Top}$_0$. Consequently, a natural question arises:
{\em what kind of  $T_0$ $L$-topological spaces can represent continuous $L$-dcpo's in a categorically isomorphic manner}.
 In this paper, we will examine the $L$-topological representations  for continuous $L$-dcpo's as well as   algebraic $L$-dcpo's.  This work can be considered as a continuation of the work \cite{YaoTFS}.

In theoretical computer science, specifically in domain theory, directed complete posets (or dcpo's, for short)  play a crucial role in most applications. Given a poset $P$, to address  the lack of some completeness of $P$, one usually considers the corresponding completion of $P$ which has the universal property \cite{Bishop}. For non-dcpo's,  Zhao  and Fan \cite{zhao-fan} introduced a new type of   directed completion of posets (called D-completion) and showed that the D-completion  possesses a universal property.  Zhang et al. \cite{Zhang-Shi-Li} studied the fuzzy directed completion of $\mathcal{Q}$-ordered sets,  generalizing the  D-completion of crisp posets. Recently, Zhang \cite{ZX-Zhang} provided a uniform approach to dealing with   $\mathcal{Z}$-directedness in  the many-valued context.
In \cite{Mislove}, Mislove generalized the notion of continuous dcpo's to  continuous posets, and obtain that the category of continuous dcpo's is  reflective in that of continuous posets.  Following this, Lawson in \cite{Lawson3} studied the round ideal completion  of continuous posets.
In this paper, we will study the directed completion of continuous $L$-ordered sets by utilizing   $L$-sobrification. The definition of directed $L$-subset  adopted  in this paper  comes from \cite{PartI} which   differs from that in \cite{ZX-Zhang,DX-Zhang-2018,Zhang-Shi-Li}. This work
 will promote  the close connections between $L$-topologies and $L$-orders.

This paper is organized as follows.  In Section 2, we list basic concepts and results about $L$-orders and $L$-topologies, and generalize the notions and results of continuous $L$-dcpo's to that of continuous $L$-ordered sets.  In Section 3, we recall some basic facts about $L$-sober spaces and provide characterizations  of $L$-sobrifications and $L$-sober spaces. In Section 4, we introduce locally super-compact spaces to  construct a categorical isomorphism between the category of the continuous $L$-dcpo's with Scott continuous maps and that of the locally super-compact $L$-sober spaces with continuous maps. In Section 5, we  establish   a categorical isomorphism  between the category of  the algebraic $L$-dcpo's with Scott continuous maps and  that of the strong locally super-compact $L$-sober spaces with continuous maps. In Section 6, we show that the $L$-sobrification of a continuous $L$-ordered set equipped  with the Scott $L$-topology is locally super-compact. Thus,  in the category of continuous  $L$-ordered sets,  we  deduce that directed completions coincide with  $L$-sobrifications.

\section{Preliminaries}

We refer to \cite{Gierz} for the content of lattice theory, to \cite{MVTop} for that of $L$-subsets and $L$-topological spaces, to \cite{Well} for that of category theory, and to \cite{PartI,YaoTFS,FanZhang} for that of $L$-ordered sets.

In this paper, $L$ always denotes a frame. A complete lattice $L$ is called a {\it frame}, or a {\it complete Heyting algebra}, if it satisfies the first infinite distributive law, that is,
$a\wedge\bigvee S=\bigvee_{s\in S}a\wedge s\ (\forall a\in L,\ \forall S\subseteq L).$
The related  adjoint pair $(\wedge,\to)$ satisfies that
$a\wedge c\leq b\Longleftrightarrow c\leq a\rightarrow b\ (\forall a,b,c\in L).$

Every map $A:X\longrightarrow L$ is called an $L$-{\it subset} of $X$, denoted by $A\in L^X$. An $L$-subset $A$ is said to be {\it nonempty} if $\bigvee_{x\in X}A(x)=1$. Let $Y\subseteq X$ and $A\in L^X$, define $A|_{Y}\in L^Y$ by $A|_{Y}(y)=A(y)$ $(\forall y\in Y)$.
For an element $a\in L$, the notation $a_X$ denotes the constant $L$-subset of $X$ with  the value $a$, i.e., $a_X(x)=a\ (\forall x\in X)$.
For all $a\in L$ and $A\in  L^{X}$,  write $a\wedge A$, for the $L$-subset given by $(a\wedge A)(x) = a\wedge A(x)$.

Let $f:X\longrightarrow Y$ be a map between two sets. The {\it Zadeh extensions} $f^\rightarrow:L^X\longrightarrow L^Y$ and $f^\leftarrow:L^Y\longrightarrow L^X$ are respectively given by
\vskip 4pt
\centerline{$f^\rightarrow(A)(y)=\bigvee\limits_{f(x)=y}A(x)\ (\forall A\in L^X)$,\   \ \qquad $f^\leftarrow(B)=B\circ f\ (\forall B\in L^Y).$}
\vskip 5pt

A subfamily $\mathcal{O}(X)\subseteq L^X$ is called an {\em $L$-topology} if
{\rm(O1)} $A,B\in\mathcal{O}(X)$ implies $A\wedge B\in \mathcal{O}(X)$;
{\rm(O2)} $\{A_j|\ j\in J\}\subseteq\mathcal{O}(X)$ implies $\bigvee_{j\in J}A_j\in\mathcal{O}(X)$;
{\rm(O3)} $a_X\in\mathcal{O}(X)$ for every $a\in L$.
The pair $(X,\mathcal{O}(X))$ is called an {\em  $L$-topological space}; elements of $\mathcal{O}(X)$ are called {\em open sets}.
Its standard name is {\em stratified $L$-topological space}. While in this paper, every $L$-topology is always assumed to be stratified, so we omit the word ``stratified''.
 As usual, we often write $X$ instead of $(X, \mathcal{O}(X))$ for an $L$-topological space.
The {\em interior operator} $(\text{-})^\circ: L^{X}\longrightarrow L^{X}$ of an $L$-topological space $X$  is defined by
$A^{\circ}=\bigvee\{B\in \mathcal{O}(X)\mid  B\leq A\}$.

A subfamily $\mathcal{B}\subseteq\mathcal{O}(X)$ is called a {\it base} of $X$ if for every $A\in\mathcal{O}(X)$, there exists $\{(B_j,a_j)|\ j\in J\}\subseteq\mathcal{B}\times L$ such that $A=\bigvee_{j\in J}a_j\wedge B_j$, or equivalently, for every $A\in\mathcal{O}(X)$,  one has $A=\bigvee_{B\in\mathcal{B}}{\rm sub}_{X}(B,A)\wedge B$ (cf. Example \ref{ex-L-ord}(3) for the operation ${\rm sub}_X$).
An $L$-topological space $X$ is called $T_0$ if $A(x)=A(y)\ (\forall A\in\mathcal{O}(X))$ implies $x=y$.
A map $f:X\longrightarrow Y$ between two $L$-topological spaces is called {\it continuous} if $f^\leftarrow(B)\in\mathcal{O}(X)$ for every $B\in\mathcal{O}(Y)$. Thus,  the continuous map $f$  gives rise to a  map from $\mathcal{O}(Y)$ to $\mathcal{O}(X)$ which sends $B$ to $f^\leftarrow(B)$. When no confuse can arise, we also use
$f^{ \leftarrow}:\mathcal{O}(Y)\longrightarrow\mathcal{O}(X)$ to denote this map.

A map $e:P\times P\longrightarrow L$ is called an $L$-{\em order} if for all $x, y, z\in P$,
(1) $e(x,x)=1$;
 (2) $e(x,y)\wedge e(y,z)\leq e(x,z)$;
and  (3) $e(x,y)\wedge e(y,x)=1$ implies $x=y$.
The pair $(P,e)$ is called an $L$-{\em ordered set}, or $L$-{\em poset}. It is customary to write $P$ for the pair $(P, e)$.
A map  $f: P\longrightarrow Q$  between two $L$-ordered sets is said to be {\em $L$-order-preserving} if for all $x, y\in P$, $e_{P}(x, y)\leq e_{Q}(f(x), f(y))$.  A map  $f: P\longrightarrow Q$  is called an {\em $L$-order-isomorphism} if $f$ is a bijection and  $e_{P}(x, y)= e_{Q}(f(x), f(y))$ for all $x, y\in P$, denoted by $P\cong Q$.
\begin{ex}\label{ex-L-ord}{\em(\cite{Lai-Zhang-2006,YaoTFS})}

{\rm(1)} Define $e_L:L\times L\longrightarrow L$ by $e_L(x, y)= x\to y$ $(\forall x, y\in L)$.   Then $e_L$ is an $L$-order on $L$.

{\rm(2)} Let $X$ be a $T_0$ $L$-topological space. Define $e_{\mathcal{O}(X)}:X\times X\longrightarrow L$ by $$e_{\mathcal{O}(X)}(x,y)=\bigwedge\limits_{A\in\mathcal{O}(X)}A(x)\rightarrow A(y)\ (\forall x,y\in X).$$ Then $e_{\mathcal{O}(X)}$ is an $L$-order on $X$, called the {\em specialization $L$-order} of $(X,\mathcal{O}(X))$. In this paper, we denote the $L$-ordered set $(X, e_{\mathcal{O}(X)})$ by $\Omega_L X$.

{\rm(3)}  Define ${\rm sub}_X :L^X\times L^X\longrightarrow L$ by
 $${\rm sub}_X(A,B)=\bigwedge\limits_{x\in X}A(x)\rightarrow B(x)\ (\forall A, B\in L^X).$$ Then ${\rm sub}_X$ is an $L$-order on $L^X$, which is called the  {\rm   inclusion $L$-order} on $L^X$.
If the background set is clear,  we always omit  the subscript in ${\rm sub}_X$.

\end{ex}

\vskip 6pt

Define ${\uparrow}A,\ {\downarrow} A\in L^P$ respectively by
\vskip 3pt
\centerline{${\uparrow}A(x)=\bigvee\limits_{y\in P}A(y)\wedge e(y,x),\ \ \ {\downarrow} A(x)=\bigvee\limits_{y\in P}A(y)\wedge e(x,y),$}
\noindent and define ${\uparrow}x$ and ${\downarrow}x$ respectively by ${\uparrow}x(y)=e(x,y)$, ${\downarrow}x(y)=e(y,x)\ (\forall x,y\in P)$. An $L$-subset $S\in L^P$ is called a {\it lower set} (resp., an {\it upper set}) if $S(x)\wedge e(y,x)\leq S(y)$ (resp., $S(x)\wedge e(x,y)\leq S(y)$) for all $x,y\in P$. Clearly, ${\downarrow}A$ and ${\downarrow}x$ (resp., ${\uparrow}A$ and ${\uparrow}x$) are lower (resp., upper) sets for all $A\in L^P$ and $x\in P$.
\vskip 6pt

Let $P$ be an $L$-ordered set. An element $x\in P$ is called a {\it supremum} of $A\in L^P$, in symbols $x=\sqcup A$, if $e(x,y)={\rm sub}(A,{\downarrow}y)\ (\forall y\in P)$.
Dually, an element $x$ is called an {\it infimum} of $A\in L^P$, in symbols $x=\sqcap A$, if $e(y,x)={\rm sub}(A,{\uparrow}y)\ (\forall y\in P)$.
If the supremum (resp., infimum) of an $L$-subset exists, then the supremum  (resp., infimum) must be unique.

 Let $P$ be an $L$-ordered set. A nonempty $L$-subset $D\in L^P$ is said to be {\it directed} if
$D(x)\wedge D(y)\leq\bigvee\limits_{z\in P}D(z)\wedge e(x,z)\wedge e(y,z)$ $(\forall x,y\in P)$.
An $L$-subset $I\in L^P$ is called an {\it ideal} of $P$ if it is a directed lower set.
Denote by $\mathcal{D}_L(P)\ ($resp., $Idl_L(P))$ the set of all directed $L$-subsets (resp., ideals) of $P$. We denote the collection of all  directed $L$-subsets (resp., ideals) with  a supremum  by $\mathcal{D}_{L}^{\ast}(P)$ (resp., $Idl_{L}^{\ast}(P))$.
An $L$-ordered set $P$ is called an $L$-{\em dcpo} if every directed $L$-subset has a supremum,
or equivalently, every ideal has a supremum.

The research on lattice-valued continuous posets (which may not be dcpo's) can be traced back to   the work of Waszkiewicz \cite{Waszkiewicz}. Inspired by Waszkiewicz's work, we adopt frame $L$ as the truth value table and extend the notions and conclusions of $L$-dcpo's (see \cite{PartI,PartII,FanZhang})  to those of  $L$-ordered sets. The primary difference between the definitions we will introduce and the earlier notions  of $L$-dcpo's is that we no longer assume the underlying $L$-ordered sets are
$L$-dcpo's. The  verification  of the results for $L$-ordered sets is analogous  to that for the corresponding results of $L$-dcpo's.

\begin{dn}
A map $f:P\longrightarrow Q$ between two $L$-ordered sets is called {\rm Scott continuous} if for every $D\in\mathcal{D}_L^{\ast}(P)$, then the supremum of $f^\rightarrow(D)$  exists and $f(\sqcup D)=\sqcup f^\rightarrow(D)$.
\end{dn}
The above definition is equivalent to that for every $I\in Idl_L^{\ast}(P)$, the supremum of $f^\rightarrow(I)$ exists and $f(\sqcup I)=\sqcup f^\rightarrow(I)$.

\begin{dn}\label{dn-con-po}
Let $P$ be an  $L$-ordered set. For every $x\in P$, define $\dda x\in L^P$ by
\begin{align*}
\dda x(y)&=\bigwedge\limits_{D\in \mathcal{D}_L^{\ast}(P)}e(x,\sqcup D)\rightarrow (\bigvee_{d\in P}D(d)\wedge e(y,d))\\
&=\bigwedge\limits_{I\in Idl_L^{\ast}(P)}e(x,\sqcup I)\rightarrow I(y).
\end{align*}
 An  $L$-ordered set $P$ is said to be {\em continuous }
if $\dda x\in\mathcal{D}_L^{\ast}(P)$ and $\sqcup{\dda}x=x$ for every $x\in P$. A continuous $L$-ordered set $P$ is called a  {\em continuous $L$-dcpo}, or an {\em $L$-domain}, if  $P$ is an $L$-dcpo.
\end{dn}

 \begin{pn}\label{lm-aux-way} Let $P$ be an $L$-ordered set. Then
$(1)$ $\forall x\in P$, $\dda x\leq {\downarrow} x$; and
$(2)$ $\forall x, u, v, y\in P$, $e(u, x)\wedge \dda y(x)\wedge e(y, v)\leq\dda v(u)$.
Thus, for every $x\in P$, $\dda x$ is a lower set.

\end{pn}
\noindent{\bf Proof.} Similar to the proof of  {(\cite[Proposition 5.7]{PartI})}.  \hfill$\Box$

 \begin{pn}\label{lm-con-int} If $P$ is a continuous $L$-ordered set, then $\dda y(x)=\bigvee_{z\in P}\dda y(z)\wedge\dda z(x)$ for all $x,y\in P$.
\end{pn}
\noindent{\bf Proof.} Similar to the proof of  {(\cite[Theorem 5.9]{PartI})}.  \hfill$\Box$

\begin{dn}
For an $L$-ordered set $P$, $A\in L^P$ is called {\em a Scott open set} if
$A$ is an upper set and $A(\sqcup D)=\bigvee_{x\in P} A(x)\wedge D(x)$ for every $D\in\mathcal{D}_L^{\ast}(P)$.
The family $\sigma_L(P)$ of all Scott open sets of $P$ forms an $L$-topology, called the {\em Scott $L$-topology} on $P$. In this paper, we use $\Sigma_L P$ to denote $(P, \sigma_L(P))$.
\end{dn}
 \begin{lm}\label{EX-sinp-dom}{\em(\cite{YaoTFS})} The $L$-ordered set $(L, e_L)$ is a continuous  $L$-dcpo. Every member in $\sigma_L(L)$ has the form $(a_{L}\wedge id_L)\vee b_L$ for a unique pair $(b, a)\in L\times L$ with $b\leq a$.
\end{lm}

\begin{pn}\label{lm-uu-op}  Let $P$ be a continuous  $L$-ordered set. Then

$(1)$ for all $x\in P$, $\uuar x$ is a Scott open set, where $\uuar x\in L^P$ is defined by $\uuar x(y)=\dda
y(x)\ (\forall y\in P)$;

$(2)$ for every $x\in P$, $({\uparrow} x)^{\circ}=\uuar x$;

$(3)$ every  Scott open set $U$ can be represented as $U=\bigvee_{x\in P}U(x)\wedge \uuar x$; thus, $\{\uuar x\mid x\in P\}$ is a base for {\color{red}$\sigma_L(P)$};

$(4)$  $(P, \sigma_L(P))$ is $T_0$.

\end{pn}
\noindent{\bf Proof.} Similar to the proof of  {(\cite[Theorems 3.7, 3.9, Proposition 3.11]{PartII})}.  \hfill$\Box$
\begin{pn}\label{pn-SO-id}
Let $P$ be a continuous $L$-ordered set. Then $\Omega_L\Sigma_L (P, e)=(P, e)$.
\end{pn}
\noindent{\bf Proof.} The proof is similar to that of \cite[Theorem 4.8]{YaoTFS}.\hfill$\Box$
\begin{pn}\label{pn-scot-map} $(1)$ Let $P$ and $Q$ be two $L$-ordered sets. If a map $f: P\to Q$ is  Scott continuous, then $f:\Sigma_L P\to \Sigma_L Q$ is  continuous.

$(2)$ Let $P$ be a continuous  $L$-ordered set and let $Q$ be a continuous  $L$-dcpo. If a map $f:\Sigma_L P\to \Sigma_L Q$ is continuous, then $f: P\to Q$ is  Scott continuous.\end{pn}

\noindent{\bf Proof.} Similar to the proof of  {(\cite[Proposition 4.5]{YaoTFS})}.  \hfill$\Box$

\section{$L$-sobriety revisited}
Since $L$-sobriety, introduced by Zhang \cite{De-Zhang}, plays a crucial role  in this paper, we recall the notions about $L$-sober spaces in this section and provide some equivalent characterizations of $L$-sobrification and $L$-sober spaces.  We first list some   basic facts about the adjunction $\mathcal{O}\dashv {\rm pt}_L$.  For notions and details not explicitly presented here,  please refer to \cite{YaoFrm, De-Zhang}.


From \cite{YaoFrm, De-Zhang}, one knows that  $L$-sobriety can be  described via an adjunction
 $$\mathcal{O}\dashv{\rm pt}_L:L\mbox{-}{\bf Top}\rightharpoonup L\mbox{-}{\bf Frm}^{op}$$
 between the category of  $L$-topological spaces and  opposite of the category of $L$-frames.
For each $L$-topological space $X$,  $(\mathcal{O}(X), {\rm sub})$ is an $L$-frame. The functor $\mathcal{O}:L\mbox{-}{\bf Top}\longrightarrow L\mbox{-}{\bf Frm}^{op}$ sends every  $L$-topological space $X$ to the   $\mathcal{O}(X)$ and sends every continuous map $f: X\longrightarrow Y$ to the $L\mbox{-}{\bf Frm}^{op}$-morphism.
$\mathcal{O}(f):\mathcal{O}(X)\longrightarrow\mathcal{O}(Y)$, where $\mathcal{O}(f)^{op}=f^{\leftarrow}: \mathcal{O}(Y)\longrightarrow\mathcal{O}(X)$.

For an $L$-topological space $X$, the functor
${\rm pt}_{L}:L\mbox{-}{\bf Frm}^{op}\longrightarrow L\mbox{-}{\bf Top}$  sends the $L$-frame $\mathcal{O}(X)$ to the $L$-topological space ${\rm pt}_{L}\mathcal{O}(X)$, where ${\rm pt}_{L}\mathcal{O}(X)$ consists of  all $L$-frame homomorphisms from $(\mathcal{O}(X), {\rm sub})$ to $(L, e_{L})$.  The related $L$-topology  on ${\rm pt}_{L}\mathcal{O}(X)$ is  $\{\phi(A)\mid A\in \mathcal{O}(X)\}$, where $\phi(A):{\rm pt}_{ L}\mathcal{O}(X)\longrightarrow L$ is defined by $p\mapsto p(A)$ for all $A\in \mathcal{O}(X)$.  This $L$-topology  is $T_0$ and is denoted by $\mathcal{O}{\rm pt}_{ L}\mathcal{O}(X)$.  In this paper, we  refer to  it as the {\it spectral $L$-topology} on ${\rm pt}_{ L}\mathcal{O}(X)$. When ${\rm pt}_{ L}\mathcal{O}(X)$ is considered as an $L$-topological space, the associated  $L$-topology is always assumed  to be the spectral $L$-topology. Therefore, we often write ${\rm pt}_{ L}\mathcal{O}(X)$ instead of $({\rm pt}_{ L}\mathcal{O}(X), \mathcal{O}{\rm pt}_{ L}\mathcal{O}(X))$, omitting the related  spectral $L$-topology.

For each  continuous map $f:X\longrightarrow Y$,  the functor  ${\rm pt}_{L}:L\mbox{-}{\bf Frm}^{op}\longrightarrow L\mbox{-}{\bf Top}$ sends the $L\mbox{-}{\bf Frm}^{op}$-morphism $\mathcal{O}(f)={f^{\leftarrow}}^{op}:\mathcal{O}(X)\longrightarrow\mathcal{O}(Y)$ to the continuous map ${\rm pt}_{ L}\mathcal{O}(f): {\rm pt}_L\mathcal{O}(X)\longrightarrow {\rm pt}_L\mathcal{O}(Y)$   defined by ${\rm pt}_{ L}\mathcal{O}(f)(p)=p\circ f^{\leftarrow}.$

 In  \cite{YaoFrm}, Yao called each member of ${\rm pt}_{L}\mathcal{O}(X)$ an
$L$-fuzzy point of $(\mathcal{O}(X), {\rm sub})$. In  this paper, we refer to   $L$-fuzzy points simply as points. The definition is restated  as follows.

\begin{dn}\label{dn-cm-fi}  A map $p:\mathcal{O}(X)\longrightarrow L$ is called a {\rm   point} of $ \mathcal{O}(X)$ if it the following conditions:

{\em (Lpt1)}  $p(A\wedge B)=p(A)\wedge p(B)$ $(\forall A, B\in \mathcal{O}(X))$;

{\em (Lpt2)}  $p(\sqcup\mathcal{A})=\sqcup p^{\rightarrow}(\mathcal{A})$ $(\forall \mathcal{A}\in L^{\mathcal{O}(X)})$.

According to {\em \cite[Proposition 5.2]{YaoFrm1}}, {\em (Lpt2)}  is equivalent to  {\em (Lpt2$^{\prime}$)} {\em +} {\em (Lpt3)}, where:


{\em (Lpt2$^{\prime}$)}  $p(\bigvee_{i\in I}A_{i})=\bigvee_{i\in I}p(A_{i})$ $(\forall \{A_{i}\mid i\in I\}\subseteq\mathcal{O}(X))$;

{\em (Lpt3)} $p(\lambda_X)=\lambda$ $(\forall \lambda\in L)$.

 The collection of all points of $\mathcal{O}(X)$ is denoted by  ${\rm pt}_{ L}\mathcal{O}(X)$.

\end{dn}
\begin{pn}\label{lm-pt} Let $X$ be an $L$-topological space.

$(1)$ For $p\in {\rm pt}_{\sss L}\mathcal{O}(X)$,  $p$ is an upper set in  $(\mathcal{O}(X), {\rm sub}_X)$;

$(2)$  $({\rm pt}_{ L}\mathcal{O}(X), {\rm sub}_{\mathcal{O}(X)})$ is an $L$-dcpo.
\end{pn}
\noindent{\bf Proof.} (1) It is clear that every point of $\mathcal{O}(X)$ is a stratified $L$-filter. By \cite[Proposition 2.16]{PartII}, $p$ is an upper set in  $(\mathcal{O}(X), {\rm sub}_X)$.

(2) Let $\mathbb{D}: {\rm pt}_{ L}\mathcal{O}(X)\longrightarrow L$ be a directed $L$-subset of ${\rm pt}_{ L}\mathcal{O}(X)$. Similar to the proof of \cite[Example 5.5] {PartI}, we have $\sqcup\mathbb{D}=\bigvee_{p\in {\rm pt}_{ L}\mathcal{O}(X)}\mathbb{D}(p)\wedge p$.
 \hfill$\Box$

In this paper,
we always assume that ${\rm pt}_{ L}\mathcal{O}(X)$, when viewed as an $L$-ordered set, is equipped with the inclusion $L$-order. For simplicity, we often write ${\rm pt}_{ L}\mathcal{O}(X)$  for $({\rm pt}_{ L}\mathcal{O}(X), {\rm sub}_{\mathcal{O}(X)})$.

\begin{ex}\label{ex-prin-fit} {\rm (1)} Let $X$ be an $L$-topological space.
For every $x\in X$, define $[x]:\mathcal{O}(X)\longrightarrow L$ by $[x](A)=A(x)$. It is easy to check  $[x]\in {\rm pt}_{ L} \mathcal{O}(X)$.

{\rm (2)} Let $X=[0, 1)$ and $L=[0,1]$. Then $\mathcal{O}(X)=\{(a_{X}\wedge id_{X})\vee b_X\mid \exists a,b\in X, b\leq a\}$ is an $L$-topology on $X$. Intuitively, every member of $\mathcal{O}(X)$ has the form shown in Figure 3.1.
Define $p: \mathcal{O}(X)\longrightarrow L$ by  $p((a_{X}\wedge id_{X})\vee b_X)=a$ for every $(a_{X}\wedge id_{X})\vee b_X\in \mathcal{O}(X)$. It is routine to check that $p$ is a point.  For $A=(1_{X}\wedge id_{X})\vee (\frac{1}{2})_X\in \mathcal{O}(X)$, we have $p(A)=1$.  However, $A(x)<1$ for every $x\in X$. This shows that  for every $x\in X$, $p\neq [x]$.
\begin{center}
  \includegraphics[height=1.5in,width=1.5in]{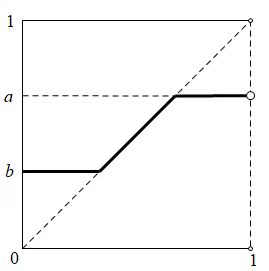}\\
 {\scriptsize {\rm {\bf Fig. 3.1:}} $(a_{X}\wedge id_{X})\vee b_X$}
 \end{center}

\end{ex}

%

 The map $\eta_{X}: X\longrightarrow {\rm pt}_{L}\mathcal{O}(X)$, defined by $x\mapsto[x]$ for all $x\in X$, is the component of the unit $\eta$ of the adjunction $\mathcal{O}\dashv{\rm pt}_L$.
 The counit  $\varepsilon$ of this adjunction gives rise to a map $\varepsilon_{\mathcal{O}(X)}^{op}:\mathcal{O}(X)\to\mathcal{O} {\rm pt}_{ L}\mathcal{O}(X)$  defined by $A\mapsto \phi(A)$ for all $A\in \mathcal{O}(X)$,
 where $\varepsilon_{\mathcal{O}(X)}^{op}$ is an isomorphism in $L\mathbf{\text{-}Frm}$. By  \cite[3.4. Idempotent pair of adjoints]{Clarka}, it is clear  that the adjunction  $\mathcal{O}\dashv {\rm pt}_L$ is idempotent.

\begin{dn}\label{dn-L-sob} {\rm (\cite{De-Zhang})} The $L$-topological space $X$ is said to be {\rm $L$-sober} if $\eta_{X}:X\longrightarrow {\rm pt}_{ L}\mathcal{O}(X)~(x\mapsto[x])$ is a bijection; hence a homeomorphism.
\end{dn}
The definition of $L$-sober spaces first appeared in \cite{De-Zhang}. These spaces are renamed  {\rm modified $L$-sober  space} in \cite{Pultr2}. For an $L$-sober space, it follows from Example \ref{ex-L-ord}(2) that
 $$\Omega_L X\cong  ({\rm pt}_{ L}\mathcal{O}(X), {\rm sub}).$$
 By Proposition \ref{lm-pt}(2), we have $\Omega_L X$ is an $L$-dcpo.


 As a special case of \cite[Proposition 5.4]{YaoFrm1}, we know that  ${\rm pt}_{L}\mathcal{O}(X)$ is  $L$-sober.  In classical  domain theory, it is known  that the Scott space of a continuous dcpo  is sober. In  the frame-valued setting, Yao obtained  the corresponding result (\cite[Proposition 5.1]{YaoFrm}).  For the ease of reference, we  list it as a  lemma below.
\begin{lm}\label{lm-dom-sob} If $P$ is a continuous $L$-dcpo, then $\Sigma_L P$  is $L$-sober.
\end{lm}
This paper will demonstrate that, just as sobriety plays an important role in domain theory, frame-valued sobriety introduced by \cite{De-Zhang} is also essential in the study of frame-valued domain theory.

 By \cite[3.6. Idempotent pairs and equivalences]{Clarka}, or  by \cite[Theorem 7.1]{Noor} one obtains  an adjunction ${\rm pt}_{L}\mathcal{O}\vdash i:L\mbox{-}{\bf Top}\rightharpoonup L\mbox{-}{\bf SobTop}$, where $i$ is the inclusion functor.  In other words, $L\mbox{-}{\bf SobTop}$ is a reflective full subcategory of $L\mbox{-}{\bf Top}$.
 The monad arising from  $\mathcal{O}\dashv {\rm pt}_L$ coincides with that arising from ${\rm pt}_{L}\mathcal{O}\vdash i$, denoted by ${\bf pt}_{L}$.It is evident  that ${\bf pt}_{L}$ is an idempotent monad (see \cite{Masahito} for idempotent monads). By \cite[Lemma 8.2]{Masahito},  the algebras of the monad ${\bf pt}_{L}$ are precisely the $L$-sober spaces.

\begin{dn}
Let $X$ be an  $L$-topological space, let $X^{S}$  be an $L$-sober space and  let $j:X\longrightarrow X^S$ be a continuous map. Then  $(X^{S}, j)$, or  $X^{S}$, is called an $L$-sobrification of $X$ if for every $L$-sober space $Y$ and every continuous map $f:X\longrightarrow Y$, there exists a unique  continuous map $\overline{f}: X^{S}\longrightarrow Y$ such that $f=\overline{f}\circ j$.

\end{dn}
 It is well-known that  $ {\rm pt}_{L}\mathcal{O}(X)$, with the map $\eta_{X}$, is an  $L$-sobrification of $(X, \mathcal{O}(X))$. By the universal property of $L$-sobrification,  $L$-sobrifications of an $L$-topological space are unique up to  homeomorphism.

The  monograph \cite{Gierz} provides  the definitions  of quasihomeomorphism and strict embedding (see \cite[Definition V-5.8]{Gierz}). By  \cite[Exercise V-5.32]{Gierz},  in the category of $T_0$ topological spaces, quasihomeomorphisms and strict embeddings are equivalent. In this paper, we extend these notions  to the category of
$L$-topological spaces.
\begin{dn}\label{dn-str-em} Let  $f:X\longrightarrow Y$ be a  map between two   $L$-topological spaces.
 Then:
 \begin{itemize}
\item $f$ is called a  {\rm quasihomeomorphism} if
$f^{\leftarrow}:\mathcal{O}(Y)\longrightarrow  \mathcal{O}(X)$ is a bijection;
\item $f: X \longrightarrow Y$  is called
a {\rm strict embedding} if it is both a quasihomeomorphism and a subspace embedding.
\end{itemize}
\end{dn}

\begin{rk}\label{rk-quasi-emb}
{\rm (1)} For a  quasihomeomorphism  $f:X\longrightarrow Y$, since $f^{\leftarrow}$ is both a  morphism in the category $L\mathbf{\text{-}Frm}$ and a bijection,  it follows  that  $f^{\leftarrow}:\mathcal{O}(Y) \longrightarrow  \mathcal{O}(X)$ is an isomorphism in   $L\mathbf{\text{-}Frm}$, equivalently, an $L$-order-isomorphism.

{\rm (2)} For every  $L$-topological space $X$, $\eta_{X}:X\longrightarrow {\rm pt}_L\mathcal{O}(X)$ is a quasihomeomorphism.
In fact,
for every $A\in \mathcal{O}(X)$,  $\eta_{X}^{\leftarrow}(\phi(A))=A$ and $\phi\circ \eta_{X}^{\leftarrow}(\phi(A))=\phi(A)$. Thus, $\eta_{X}^{\leftarrow}\circ \phi=id_{\mathcal{O}(X)}$ and $\phi\circ\eta_{X}^{\leftarrow}=id_{\mathcal{O}pt\mathcal{O}(X)}$.

{\rm (3)} Let $f:X\longrightarrow Y$ be a  map between two $T_0$ $L$-topological spaces. It is routine to verify that   $f$ is a subspace embedding if and only if $f^{\leftarrow}:\mathcal{O}(Y)\longrightarrow  \mathcal{O}(X)$ is a surjection. Therefore, in the category of $T_0$ $L$-topological spaces,  quasihomeomorphisms are equivalent to  strict embeddings.

\end{rk}

In what follows, we give    characterizations of $L$-sobrification.

\begin{tm}\label{tm-char-sobri}
Let $X$ be an $L$-topological space and let $Y$ be an $L$-sober space. Then the following  are equivalent:

$(1)$ $Y$ is an $L$-sobrification of $X$;

$(2)$ there exists a quasihomeomorphism $j:X \longrightarrow Y$;

$(3)$ $(\mathcal{O}(X), {\rm sub}_X)$  and $(\mathcal{O}(Y), {\rm sub}_Y)$ are $L$-order-isomorphic.
\end{tm}
\noindent {\bf Proof.}
$(1)\Rightarrow (2)$: Since  $Y$ and  ${\rm pt}_L\mathcal{O}(X)$ are  $L$-sobrifications of $X$, we know that ${\rm pt}_L\mathcal{O}(X)$ is homeomorphism to $Y$. Notice that $\eta_{X}:X\longrightarrow {\rm pt}_L\mathcal{O}(X)$ is a quasihomeomorphism. Hence, there exists a quasihomeomorphism $j:X \longrightarrow Y$.

$(2)\Rightarrow (3)$: Clear from Remark \ref{rk-quasi-emb}(1).

$(3)\Rightarrow(1)$: By Definition \ref {dn-L-sob} and the constructions of $\mathcal{O}{\rm pt}_L\mathcal{O}(X)$ and $\mathcal{O}{\rm pt}_L\mathcal{O}(Y)$, or by the fact that ${\rm pt}_{L}:L\mbox{-}{\bf Frm}^{op}\longrightarrow L\mbox{-}{\bf Top}$ is a functor and $(\mathcal{O}(X), {\rm sub}_X)$ is isomorphic to $(\mathcal{O}(Y), {\rm sub}_Y)$ in the category $L\mathbf{\text{-}Frm}^{op}$, we conclude that
$$({\rm pt}_L\mathcal{O}(X), \mathcal{O}{\rm pt}_L\mathcal{O}(X))\cong({\rm pt}_L\mathcal{O}(Y), \mathcal{O}{\rm pt}_L\mathcal{O}(Y)).$$
Since $Y$ be an $L$-sober space, we have $(Y, \mathcal{O}(Y))\cong({\rm pt}_L\mathcal{O}(Y), \mathcal{O}{\rm pt}_L\mathcal{O}(Y))$.
Thus, $(Y, \mathcal{O}(Y))\cong ({\rm pt}_L\mathcal{O}(X), \mathcal{O}{\rm pt}_L\mathcal{O}(X))$. Notice that    ${\rm pt}_L\mathcal{O}(X)$ is an   $L$-sobrification of $X$. Hence,  $Y$ is an $L$-sobrification of $X$.\hfill$\Box$

In the classical setting, Lawson  typically took  Theorem \ref{tm-char-sobri}(2) as the definition of  sobrification (see \cite{Lawson2}).
The following example, along with the example in Section 6, demonstrates  that   Theorem \ref{tm-char-sobri}(3) can be effectively used to identify the  $L$-sobrification of an $L$-topological space.
\begin{ex} Let $X$, $L$  and  $\mathcal{O}(X)$ be the same as in Example \ref{ex-prin-fit}(2). By Lemma \ref{EX-sinp-dom}, it is easy to see that $(\mathcal{O}(X), {\rm sub}_X)\cong (\sigma_L(L), {\rm sub}_L)$. It follows from Lemma \ref{EX-sinp-dom} and Lemma \ref {lm-dom-sob} that $(L, \sigma_L(L))$ is an $L$-sober space. So,  by Theorem  \ref{tm-char-sobri}, $(L, \sigma_L(L))$ is an $L$-sobrification of $(X, \mathcal{O}(X))$.
\end{ex}

\begin{pn}\label{pn-flat-homo} Let $X$, $Y$ be two $L$-topological spaces and let $f:X\longrightarrow Y$ be a continuous map. Then following  are equivalent:

$(1)$  $f$ is a quasihomeomorphism;

$(2)$   ${\rm pt}_L\mathcal{O}(f): {\rm pt}_L\mathcal{O}(X)\longrightarrow {\rm pt}_L\mathcal{O}(Y)$ is a homeomorphism.

\end{pn}
\noindent {\bf Proof.} $(1)\Rightarrow (2)$: It is clear since $\mathcal{O}(f)$ is an isomorphism   in $L\mathbf{\text{-}Frm}^{op}$.

$(2)\Rightarrow (1)$:
Since ${\rm pt}_L\mathcal{O}(f)$ is a homeomorphism, it follows that $\mathcal{O}{\rm pt}_L\mathcal{O}(f): \mathcal{O}{\rm pt}_L\mathcal{O}(X)\longrightarrow \mathcal{O}{\rm pt}_L\mathcal{O}(Y)$ is an isomorphism in $L\mathbf{\text{-}Frm}^{op}$. In $L\mathbf{\text{-}Frm}$, the adjunction $\mathcal{O}\dashv {\rm pt}_L$ gives rise to the following diagram.
\begin{displaymath}
\xymatrix@=8ex{\mathcal{O}(Y)\ar[r]^{\varepsilon_{\mathcal{O}(Y)}^{op}\ \ }_{\cong\ \ \ }\ar[d]_{\mathcal{O}(f)^{op} }&\mathcal{O}{\rm pt}_L\mathcal{O}(Y)\ar[d]^{(\mathcal{O}{\rm pt}_L\mathcal{O}(f))^{op}}_{\cong}\\
\mathcal{O}(X)\ar[r]_{\varepsilon_{\mathcal{O}(X)}^{op}\ \ }^{\cong\ \ \ }&  \mathcal{O}{\rm pt}_L\mathcal{O}(X)}
\end{displaymath}
 Therefore, $\mathcal{O}(f)^{op}=({\varepsilon_{\mathcal{O}(X)}^{op}})^{-1}\circ(\mathcal{O}{\rm pt}_L\mathcal{O}(f))^{op}\circ\varepsilon_{\mathcal{O}(Y)}^{op}$. Hence $\mathcal{O}(f)^{op}=f^{\leftarrow}$ is an isomorphism in $L\mathbf{\text{-}Frm}$. It follows that $f$ is a quasihomeomorphism.
\hfill$\Box$

 Readers can refer to
 \cite[Section II-3]{Gierz} for the notion of {\rm $J$-injective} object in terms of  category. We use  $L$-{\bf Top}$_0$ to denote the category of  $T_0$ $L$-topological spaces with continuous maps.
In $L$-$\mathbf{Top}_0$, let $J_1$  and $J_2$ be the class of all subspace embeddings and  all strict embeddings respectively.  Then a $T_0$ $L$-topological space $X$ is said to be {\it  injective} (resp., {\it strictly injective}) if  $X$  is $J_1$-injective (resp.,   $J_2$-injective) in the category $L$-$\mathbf{Top}_0$.
In the  following, we will show that $L$-sober spaces are equivalent to  strictly injective $T_0$ $L$-topological spaces. It is a lattice-valued type of  Escard\'{o} and  Flagg's corresponding result (see \cite{M.H.Es}). But Escard\'{o} and  Flagg did not give a direct proof.

\begin{tm}\label{tm-inj-sob} The $L$-topological space  $X$ is a strictly injective $T_0$ $L$-topological space if and only if $X$ is an $L$-sober space.
\end{tm}
\noindent{\bf Proof.}
{\bf Necessity.} It follows from Remark \ref{rk-quasi-emb}(2)(3) that  $\eta_{X}$  is a strict embedding.  Since $X$ is strictly injective,  there exists a continuous map $r: {\rm pt}_L\mathcal{O}(X)\longrightarrow X$ such that $r\circ \eta_X=id_X$.
By \cite[Proposition 5.8]{YaoFrm}, $X$ is $L$-sober.

{\bf Sufficiency.} Let $X$ be an $L$-sober space and     let $j: Y\longrightarrow Z$ be a strict embedding between two  $T_0$ $L$-topological spaces. Then $\eta_{X}:X\longrightarrow {\rm pt}_L\mathcal{O}(X)$ is a bijection.  By Proposition \ref{pn-flat-homo},  ${\rm pt}_L\mathcal{O}(j): {\rm pt}_L\mathcal{O}(Y)\longrightarrow {\rm pt}_L\mathcal{O}(Z)$ is a homeomorphism.  For every   continuous map $f: Y\longrightarrow X$,
construct $\overline{f}:Z\longrightarrow X$ as follows:
$$\overline{f}:Z \xrightarrow{\eta_{Z}} {\rm pt}_L\mathcal{O}(Z)\xrightarrow{({\rm pt}_L\mathcal{O}(j))^{-1}} {\rm pt}_L\mathcal{O}(Y)\xrightarrow{{\rm pt}_L\mathcal{O}(f)}{\rm pt}_L\mathcal{O}(X)\xrightarrow{(\eta_{X})^{-1}}X.$$
By the naturality of $\eta$,  one has
\begin{align*}
\overline{f}\circ j&=(\eta_{X})^{-1}\circ {\rm pt}_L\mathcal{O}(f)\circ ({\rm pt}_L\mathcal{O}(j))^{-1}\circ \eta_{Z}\circ j\\
&=(\eta_{X})^{-1}\circ {\rm pt}_L\mathcal{O}(f)\circ  ({\rm pt}_L\mathcal{O}(j))^{-1}\circ {\rm pt}_L\mathcal{O}(j)\circ \eta_{Y}\\
&=(\eta_{X})^{-1}\circ {\rm pt}_L\mathcal{O}(f)\circ \eta_{Y}\\
&=(\eta_{X})^{-1}\circ \eta_{X}\circ f=f.
\end{align*}
So, $\overline{f}\circ j=f$.
To sum up, $X$ is strictly injective.\hfill$\Box$

 By the naturality of $\eta$,  it is straightforward to verify that the  map $\overline{f}$ constructed in the proof above is unique.

\section{Topological representations of continuous $L$-dcpo's }
In \cite{YaoTFS}, Yao demonstrated that  every continuous $L$-lattice (referred to as a fuzzy
continuous lattice in \cite{YaoTFS}) equipped with Scott $L$-topology is an injective $T_0$ $L$-topological space, and conversely, the specialization  $L$-ordered set of an injective $T_0$ $L$-topological space is a continuous $L$-lattice. These transformations form a categorical isomorphism between the category of continuous $L$-lattices and that  of injective $T_0$ $L$-topological spaces. Thus one has the following corresponding classes of $L$-ordered sets and  $T_0$ $L$-topological spaces:
\vskip 5pt
 (4.1)\qquad  $\{ \mbox{continuous } L\mbox{-lattices}\}\Longleftrightarrow \{ \mbox{injective } T_0 \ L\mbox{-topological spaces}\}$
\vskip 5pt
It is thus natural to wonder, {\em  what kind of  $T_0$ $L$-topological spaces can continuous $L$-dcpo's be categorically isomorphic to?} To address this issue, we draw inspiration from the topological representation of classical continuous dcpo's. It follows from the results of \cite{Banaschewski,Erne,R. Hofmann,Lawson,XU-mao-form,Zhao-Parial} that continuous dcpo's can be represented by the following equivalent types of topological spaces, despite their differing formulations:
\begin{itemize}
 \item $X$ is  both a sober space and  also a  C-space;
 \item $X$ is   a sober space which has  an injective hull;
  \item $X$ is   a sober space whose  open set lattice is completely distributive;
  \item $X$ is  both    sober  and  locally super-compact.
\end{itemize}
For the notions of the first three spaces, refer to \cite{Zhao-Parial}; and for the last one, refer to \cite{XU-mao-form}.
In \cite{XU-mao-form}, Xu and Mao provided  representations for various  domains via the so-called formal points related to  super-compact quasi-bases of  sober spaces. This representation of  continuous dcpo's is both simple and elegant, as it is obtained through set-theoretic sequences of proof without  using  higher-order definitions or properties, such as complete distributivity of open set lattice and the categorical properties of  spaces.   In this section, we  replace  the  so-called formal points defined in \cite{XU-mao-form} with the points (equivalently, completely prime filters) of the open set  lattice. This change will make the corresponding proofs more direct and will facilitate a smooth generalization of the results to the frame-valued setting.  Xu and Mao  did not examine the relationship between various domains  and locally super-compact sober spaces from a categorical  perspective. However, in this section,  we  will establish an isomorphism between continuous $L$-dcpo's and locally super-compact sober spaces.
%

We first introduce the notions of  super-compact $L$-subsets and locally super-compact $L$-topological spaces.

\begin{dn}\label{dn-l-lsc} Let $X$ be an $L$-topological space.
 \begin{itemize}
\item An $L$-subset $A\in L^{X}$ is called a  {\rm super-compact $L$-subset}, if $\bigvee_{x\in X}A(x)=1$, and  for every family $ \{ V_i\mid i\in I\}\subseteq \mathcal{O}(X)$,
$${\rm sub}(A, \bigvee_{i\in I}V_{i})=\bigvee_{i\in I}{\rm sub}(A, V_{i}).$$
Denote all the super-compact $L$-subsets of $(X, \mathcal{O}(X))$  by ${\rm SC}(X)$.

\item The space  $X$ is called {\rm locally super-compact}, if for every $A\in \mathcal{O}(X)$,
$$A=\bigvee\limits_{B\in {\rm SC}(X)}{\rm sub}(B, A)\wedge B^{\circ}.$$
It is clear that  the family $\{B^{\circ}\mid B\in {\rm SC}(X)\}$ is a base of the locally super-compact $L$-topological space  $X$.
 \end{itemize}
\end{dn}

For $A\in L^X$, define $[A]:\mathcal{O}(X)\longrightarrow L$ by
$$[A](B)={\rm sub}(A,B)\ (\forall B\in\mathcal{O}(X)).$$
It is clear that  for each $C\in {\rm SC}(X)$,  $[C]\in {\rm pt}_{L}\mathcal{O}(X)$.

The following shows that  every   continuous $L$-ordered set can induce a locally super-compact  $L$-topological space.

\begin{pn}\label{tm-dom-sob} Let $P$  is a continuous $L$-ordered set. Then $\Sigma_L P$ is a locally super-compact  space.     Thus,   if $P$  is also an $L$-dcpo,  then $\Sigma_L P$  is  a locally super-compact  $L$-sober space.
\end{pn}
\noindent{\bf Proof.}
 For every $y\in P$, we have $\bigvee_{x\in P}{\uparrow} y(x)\geq {\uparrow} y(y)=1$.  Since every Scott open set is an upper set, we have
$$\mbox{sub}({\uparrow} y, \bigvee_{i\in I}V_{i})= \bigvee_{i\in I}V_{i}(y)= \bigvee_{i\in I}\mbox{sub}({\uparrow} y, V_{i}).$$
Therefore  ${\uparrow} y$ is super-compact.
For every $A\in \sigma_L(P)$, by Proposition  \ref{lm-uu-op}, we have
$$A=\bigvee_{y\in P} A(y)\wedge\uuar y=\bigvee_{y\in P} {\rm sub}({\uparrow} y, A)\wedge({\uparrow} y)^{\circ}.$$
 Consequently $(P, \sigma_L(P))$ is a locally super-compact  space.\hfill$\Box$

The following example indicates that, in an $L$-topological space,  there can exist  nonempty $L$-subsets that are not locally super-compact.
\begin{ex} Let $L=\{0, a, b, c, 1\}$ be a lattice with $0\leq a\leq c\leq 1$, $0\leq b\leq c\leq 1$ and $a\parallel b$. Clearly, $L$ is a frame. For $A=id_L\vee c_L$, $B=id_L\vee a_L$ and $C=id_L\vee b_L$, by   Lemma \ref{EX-sinp-dom}, $A, B, C\in \sigma_L((L, e_L))$.
Clearly,  $\bigvee_{x\in L}A(x)=1$, but
$${\rm sub}( A, B)\vee {\rm sub}( A, C)=a\vee b=c;$$
and
$${\rm sub}( A, B\vee C)={\rm sub}( A, A)=1.$$
It shows that $ A$ is not a locally super-compact $L$-subset in $(L, \sigma_L(L))$.
\end{ex}

\begin{lm}{\em (\cite[Proposition 4.3]{LiQG1})}\label{lm-base} Let $P$ be an $L$-dcpo and $x\in P$. If there exists a directed $L$-subset $D\in L^P$ such that $D\leq \dda x$ and $\sqcup D=x$, then $\dda x$ is  directed and $x=\sqcup\dda x$.
\end{lm}

\begin{pn}\label{pn-lsc-dom} Let $X$ be a  locally super-compact $L$-topological space. Then in  $({\rm pt}_{ L}\mathcal{O}(X), {\rm sub}_{\mathcal{O}(X)})$, $\dda [x]$ is directed and $[x]=\sqcup\dda [x]$ for every $x\in X$.
\end{pn}
\noindent{\bf Proof.}
By Lemma \ref{lm-base}, we only need to show that for every $x\in X$, there exists a directed $L$-subset $\mathbb{D}_{x}$ of ${\rm pt}_{ L}\mathcal{O}(X)$ such that $\mathbb{D}_{x}\leq \dda [x]$ and $\sqcup\mathbb{D}_{x}=[x]$.
For every $x\in X$, define $\mathbb{D}_{x}: {\rm pt}_{ L}\mathcal{O}(X)\longrightarrow L$ by
$$\mathbb{D}_{x}(p)=\bigvee\{A^{\circ}(x)\mid A\in {\rm SC}(X),  p=[A] \}\ (\forall p\in {\rm pt}_{ L}\mathcal{O}(X)).$$
 It is easy to see $\bigvee_{p\in {\rm pt}_{ L}\mathcal{O}(X))}\mathbb{D}_x(p)=\bigvee_{A\in {\rm SC}(X)}A^{\circ}(x)=1$.
Let $p_{1}$, $p_2\in {\rm pt}_{ L}\mathcal{O}(X)$.  Without loss of generality, we assume that there exist $A, B\in {\rm SC}(X)$ such that  $p_1=[A]$ and    $p_2=[B]$.    Then we have
\begin{align*}
&\qquad  \mathbb{D}_{x}(p_1)\wedge \mathbb{D}_{x}(p_2)\\
&= \bigvee_{\substack {C_{1}\in {\rm SC}(X)\\p_1=[C_{1}]}}C_{1}^{\circ}(x)\wedge \bigvee_{\substack{C_{2}\in {\rm SC}(X)\\ p_2=[C_{2}]}}C_{2}^{\circ}(x)\\
&= \bigvee\{C_{1}^{\circ}(x)\wedge C_{2}^{\circ}(x)\mid  C_{1},C_{2}\in {\rm SC}(X),  p_1=[C_{1}], p_2=[C_{2}]\}\\
&= \bigvee\{\bigvee_{C\in {\rm SC}(X)}C^{\circ}(x)\wedge {\rm sub}(C, C_{1}^{\circ}\wedge C_{2}^{\circ})\mid C_{1},C_{2}\in {\rm SC}(X), p_1=[C_{1}], p_2=[C_{2}]\}\\
&\leq\bigvee\{\bigvee_{C\in {\rm SC}(X)}C^{\circ}(x)\wedge {\rm sub}(C, C_{1})\wedge {\rm sub}(C,C_{2})\mid C_{1},C_{2}\in {\rm SC}(X), p_1=[C_{1}], p_2=[C_{2}]\}\\
&\leq\bigvee_{C\in {\rm SC}(X)} \mathbb{D}_{x}([C])\wedge {\rm sub}(p_1, [C])\wedge {\rm sub}(p_2,[C])\\
&=\bigvee_{p\in {\rm pt}_{ L}\mathcal{O}(X)}\mathbb{D}_{x}(p)\wedge {\rm sub}(p_1, p)\wedge {\rm sub}(p_2,p).
\end{align*}
This shows that   $\mathbb{D}_{x}$ is directed.

Now we  show that $\mathbb{D}_{x}\leq \dda [x]$. Let  $p\in {\rm pt}_{ L}\mathcal{O}(X)$.  Without loss of generality, we assume that there exists $B\in {\rm SC}(X)$ such that  $p=[B]$.
 For every ideal $\mathbb{I}$ of ${\rm pt}_{ L}\mathcal{O}(X)$ and $B\in {\rm SC}(X)$ with $p=[B]$. By Proposition \ref{lm-pt}, we have
\begin{align*}
{\rm sub}([x], \sqcup\mathbb{I})&\leq [x](B^{\circ})\to \sqcup\mathbb{I}(B^{\circ})\\
&=B^{\circ}(x)\to \bigvee_{q\in {\rm pt}_{ L}\mathcal{O}(X)}\mathbb{I}(q)\wedge q(B^{\circ})\\
&\leq B^{\circ}(x)\to \bigvee_{q\in {\rm pt}_{ L}\mathcal{O}(X)}\mathbb{I}(q)\wedge {\rm sub}([B], q)\\
&\leq B^{\circ}(x)\to \mathbb{I}([B])\\
&=B^{\circ}(x)\to \mathbb{I}(p).
\end{align*}
\noindent Therefore, $$B^{\circ}(x)\leq {\rm sub}([x], \sqcup\mathbb{I})\to \mathbb{I}(p).$$ By the arbitrariness of $\mathbb{I}$ and $B$, we have
\begin{align*}
\mathbb{D}_{x}(p)&=\bigvee\{B^{\circ}(x)\mid B\in {\rm SC}(X), p=[B]\}\\
&\leq\bigwedge\limits_{\mathbb{I}\in Idl_{L}({\rm pt}_{ L}\mathcal{O}(X))} {\rm sub}([x], \sqcup\mathbb{I})\to \mathbb{I}(p)\\
&=\dda [x](p).
\end{align*}

\noindent To the end, we need  show that $\sqcup\mathbb{D}_{x}=[x]$. For every $A\in \mathcal{O}(X)$,
\begin{align*}
\sqcup\mathbb{D}_x(A)&=\bigvee_{p\in {\rm pt}_{ L}\mathcal{O}(X)}\mathbb{D}_x(p)\wedge p(A)\\
&=\bigvee_{p\in {\rm pt}_{ L}\mathcal{O}(X)}\bigvee\{B^{\circ}(x)\wedge [B](A)\mid B\in {\rm SC}(X), p=[B]\}\\
&=\bigvee_{B\in {\rm SC}(X)}B^{\circ}(x)\wedge {\rm sub}(B, A)\\
&=A(x)
=[x](A),
\end{align*}\label{lm-base-way}
as desired. \hfill$\Box$

The following demonstrates that by adding  the condition of $L$-sobriety to a locally super-compact $L$-sober space, this space can induce a  continuous $L$-dcpo.

\begin{pn}\label{pn-sob-cont}   Let $X$ be a locally super-compact L-sober space. Then
 $\Omega_L X$ is a  continuous $L$-dcpo.

\end{pn}
\noindent{\bf Proof.} It follows from the $L$-sobriety  of $X$ that ${\rm pt}_{ L}\mathcal{O}(X)=\{[x]\mid x\in X\}$ is an $L$-dcpo. Then by Proposition \ref{pn-lsc-dom}, ${\rm pt}_{ L}\mathcal{O}(X)$ is a continuous $L$-dcpo.
Since $\Omega_L X\cong  ({\rm pt}_{ L}\mathcal{O}(X), {\rm sub})$,
 $\Omega_L X$ is a  continuous $L$-dcpo.\hfill$\Box$

\begin{tm}\label{tm-SO-id} {\rm (1)} Let $X$ be a locally super-compact $L$-sober space. Then $\Sigma_L\Omega_L X=(X, \mathcal{O}(X))$.

{\rm (2)} Let $P$ be a continuous $L$-dcpo. Then $\Omega_L\Sigma_L (P, e)=(P, e)$.
\end{tm}
\noindent{\bf Proof.} (1) It follows from  \cite[Proposition 5.5]{YaoFrm} that $ \mathcal{O}(X)\subseteq\sigma_L(\Omega_{L} X)$. In the following, we only need to show that $\sigma_{L}(\Omega_{L} X)\subseteq \mathcal{O}(X)$.

Let  $x\in X$.  Then $\sqcup\mathbb{D}_x=[x]$, where  $\mathbb{D}_x$ is the one constructed in the proof of Proposition \ref{pn-lsc-dom}.  Define $D_x\in L^{X}$ by $D_x(a)=\mathbb{D}_x([a])$ $(\forall a\in X)$. By the $L$-sobriety of $X$, $D_x$ is well-defined. Since $\Omega_{L} X\cong({\rm pt}_{ L}\mathcal{O}(X), {\rm sub})$ and $\mathbb{D}_x$ is a directed $L$-subset of $({\rm pt}_{ L}\mathcal{O}(X), {\rm sub})$, it is clear that $D_x$ is a directed $L$-subset of $\Omega_{L} X$ and $\sqcup D_x=x$. Let $A\in \sigma_{L}(\Omega_{L} X)$. Then
\begin{align*}
A(x)&=A(\sqcup D_x)\\
&=\bigvee_{a\in X}A(a)\wedge D_x(a)\\
&=\bigvee_{a\in X}[a](A)\wedge\mathbb{D}_x([a])\\
&=\bigvee_{a\in X}[a](A)\wedge\bigvee_{\substack{B\in {\rm SC}(X)\\ [a]=[B]}} B^{\circ}(x)\\
&=\bigvee_{ B\in {\rm SC}(X)} [B](A)\wedge B^{\circ}(x)\\
&=\bigvee_{ B\in {\rm SC}(X)}{\rm sub}(B, A)\wedge B^{\circ}(x),
\end{align*}
This shows that
$$A=\bigvee_{ B\in {\rm SC}(X)} {\rm sub}(B, A)\wedge B^{\circ}\in \mathcal{O}(X).$$ Thus, $\sigma_{L}(\Omega_{L} X)\subseteq \mathcal{O}(X)$, as desired.

(2) It is clear from Proposition \ref{pn-SO-id}.
\hfill$\Box$

Let $L\mathbf{\text{-}LSCSob}$ denote the category of  locally super-compact $L$-sober spaces and continuous maps, and let $L\mathbf{\text{-}CDom}$ denote the category of  continuous $L$-dcpo's and Scott continuous maps.
In the following theorem demonstrates that locally super-compact $L$-sober spaces and  continuous $L$-dcpo's are the same thing up to  a categorical isomorphism.

\begin{tm}\label{tm-iso-sob-dom}   The categories $L$-{\bf LSCSob}  and  $L$-{\bf CDom} are isomorphic.
\end{tm}
\noindent{\bf Proof.}
Define assignment $\Sigma_L:L\mathbf{\text{-}CDom}   \to L\mathbf{\text{-}LSCSob}$ by
$$(f: P\to Q)\mapsto (f: \Sigma_L P\to  \Sigma_L Q).$$
It follows from  Proposition \ref{tm-dom-sob} and Proposition \ref{pn-scot-map}(1) that $\Sigma_L$ is a functor.
Define assignment $\Omega_L: L\mathbf{\text{-}LSCSob}   \longrightarrow L\mathbf{\text{-}CDom}$ by
$$(f:X\to Y)\mapsto(f: \Omega_L X\to \Omega_L Y).$$
Let $X$, $Y$ be two locally super-compact $L$-sober spaces and let $f: X\longrightarrow Y$ be a  continuous map.  By Proposition \ref{pn-sob-cont}, $\Omega_L X$ and  $\Omega_L Y$ are continuous $L$-dcpo's.  It follows from Theorem \ref{tm-SO-id}(1) and  Proposition \ref{pn-scot-map}(2) that $f:\Omega_L
X\longrightarrow \Omega_L Y$ is Scott continuous. Thus, $\Omega_L$ is a functor. By Theorem \ref{tm-SO-id}, it is clear that $\Sigma_L\circ\Omega_L=id_{L\mathbf{\text{-}LSCSob}}$ and $\Omega_L\circ\Sigma_L=id_{L\mathbf{\text{-}CDom}}$. Thus, $L\mathbf{\text{-}LSCSob}$  and  $L\mathbf{\text{-}CDom}$ are isomorphic.
\hfill$\Box$\\

\section{Topological representations of algebraic $L$-dcpo's}

 In this section,  we will   explore a specific class  of locally super-compact  $L$-topological spaces to give representations for  algebraic $L$-dcpo's. We begin by recalling some concepts related to $L$-dcpo's, as defined in \cite{Su-Li-Intel}, which differ  from the lattice-valued type of algebraic dcpo's  introduced in \cite{LiQG1}.
\begin{dn}{\rm (\cite{Su-Li-Intel})}\label{dn-alg-dom}
 Let $P$ be an  $L$-dcpo.
  \begin{itemize}
\item The element $x\in P$ is said to be  {\rm compact} if  $\dda x(x)=1$. The set of all compact elements of $P$ is denoted by $K(P)$.
\item Define $k(x)\in L^{P}$ by $k(x)(y)=e(y, x)$ for $y\in K(P)$ and
otherwise $0$.
If for every $x\in P$, $k(x)$ is directed and $x=\sqcup k(x)$, then $P$ is called an {\rm algebraic $L$-dcpo}.
 \end{itemize}
\end{dn}

\begin{rk}\label{rk-yao-alg}
{\rm (1)} $\dda x(x)=1$ iff $I(x)=e(x, \sqcup I)$ $(\forall I\in Idl_{L}(P))$ iff

$$e(x, \sqcup D)=\bigvee_{d\in P}D(d)\wedge e(x, d)\ (\forall D\in \mathcal{D}_{L}(P)).$$


{\rm (2)} Clearly,  for every algebraic $L$-dcpo $P$, $k(x)\leq \dda x$. It follows directly from Lemma \ref{lm-base} that every  algebraic $L$-dcpo $P$ is  a continuous $L$-dcpo. However, we will show that a continuous $L$-dcpo is not necessarily an algebraic $L$-dcpo.

\end{rk}
\begin{ex}
{\rm (1)} Let  $L=[0,1]$. Then $(L, e_L)$ is a continuous $L$-dcpo. Then by the proof of {\em \cite[Theorem 3.5]{YaoTFS}}, we have $\uuar x= id_L$  for every $x\in (0, 1]$ and $\uuar 0=1_L$. Thus for every $a\in [0, 1]$, we have
$$
\ \ \dda a(x)=\left\{\begin{array}{ll}1,& x=0;\\
 a,& x\in(0,1].
\end{array}\right.
$$
It follows that  $K(L)=\{0, 1\}$. Thus
$$
\ \ k(a)(x)=\left\{\begin{array}{ll}1,& x=0;\\
0, &x\in(0,1);\\
 a,& x=1.
\end{array}\right.
$$
For every $y\in [0, 1]$, we have
\vskip 3pt
\centerline{${\rm sub}(k(a), {\downarrow} y)=(k(a)(0)\rightarrow e_L(0, y))\wedge(k(a)(1)\rightarrow e_L(1, y))=e_L(a, y).$}
\vskip 3pt
\noindent This shows that $\sqcup k(a)=a$. Thus  $(L, e_L)$ is an algebraic $L$-dcpo.

{\rm (2)} Let $X=[0, 1]\cup\{\bot\}$ and  $L=\{0, \frac{1}{2}, 1\}$. Define an $L$-order $e$ on $X$ as follows:

$$
\ \  e(x, y)=\left\{\begin{array}{ll}1,& x, y\in [0,1], x\leq y;\\
 1, & x=\bot,  y\in X;\\
\frac{1}{2}& x, y\in [0,1], y<x;\\
0 & y=\bot,  x\in[0,1].
\end{array}\right.
$$
\vskip 3pt
 It is routine to check that every  $I\in Idl_L(X)$  can be precisely described by one of the following three types of $L$-subsets, where $a\in [0,1]$.
\vskip 8pt
\centerline{$
\ I_1(x)=\left\{\begin{array}{ll}0, &x\in[0,1];\\
1,& x=\bot,
\end{array}\right.$
$
 \ I_2(x)=\left\{\begin{array}{ll}\frac{1}{2}, &x\in [a, 1];\\
1,& \mbox{otherwise},
\end{array}\right.
$
$
\ I_3(x)=\left\{\begin{array}{ll}\frac{1}{2}, &x\in (a, 1];\\
1, &\mbox{otherwise}.
\end{array}\right.
$
}

\vskip 8pt
\noindent It is easy to verify that $\sqcup I_1=\bot$ and $\sqcup I_2=\sqcup I_3=a$.

For $\bot\in X$ and for every $a\in [0,1]$, we have
\vskip 8pt
\centerline{$
\ \ \dda \bot(x)=\left\{\begin{array}{ll}0, &x\in[0,1];\\
1,& x=\bot,
\end{array}\right.
$
$
\ \ \dda a(x)=\left\{\begin{array}{ll}\frac{1}{2}&x\in[a,1];\\
1,& \mbox{otherwise}.
\end{array}\right.
$}
\vskip 8pt
\noindent Thus, $(X, e)$ is a continuous $L$-dcpo with $K(X)=\{\bot\}$. For every $a\in [0, 1]$, $k(a)=I_1$. Since $\sqcup k(a)=\bot\neq a$, it follows that $(X, e)$ is not an algebraic $L$-dcpo.

\end{ex}

In \cite{XU-mao-form}, Xu and Mao used topological   spaces with a base consisting of super-compact open sets to provide representations for algebraic dcpo's. This approach inspires us to  introduce the notion of strong locally super-compact $L$-topological spaces, defined as follows.

\begin{dn}Let $X$ be an $L$-topological space. If $X$    has   a base $\mathcal{B}$ consisting of super-compact open sets;  that is, $\mathcal{B}\subseteq {\rm SC}(X)$, then  $X$ is said to be  {\rm strong locally super-compact}.
 \end{dn}

\begin{pn} A strong locally super-compact $L$-topological space $X$ is a locally super-compact $L$-topological space.
\end{pn}
\noindent{\bf Proof.} Assume $\mathcal{B}$ is  a base of $X$ which  consists of super-compact open sets. For every $A\in \mathcal{O}(X)$, we have
\vskip 6pt
\centerline{$A=\bigvee_{B\in \mathcal{B}} {\rm sub}(B, A)\wedge B\leq \bigvee_{B\in {\rm SC}(X)}{\rm sub}(B, A)\wedge B^{\circ}\leq A.$}
Therefore $A= \bigvee_{B\in {\rm SC}(X)}{\rm sub}(B, A)\wedge B^{\circ}$.
Thus, $X$ is locally super-compact.\hfill$\Box$

The following shows that every algebraic $L$-dcpo can induce a strong locally super-compact $L$-topological space.
\begin{pn}\label{pn-alg-bas} Let $P$ be an algebraic $L$-dcpo.  Then $\Sigma_L P$ is a strong locally super-compact $L$-topological space.
\end{pn}
\noindent{\bf Proof.} We first show that  $\{{\uparrow} y\mid y\in K(P)\}$ is  a base of $\Sigma_L P$. Let $y\in K(P)$. Then for every $D\in \mathcal{D}_{L}(P)$, we have
$${\uparrow} y(\sqcup D)=e(y, \sqcup D)= \bigvee_{d\in P}e(y, d)\wedge D(d).$$
Thus ${\uparrow} y\in \sigma_{L}(P)$. Since $P$ is an algebraic $L$-dcpo,  for every $A\in \sigma_{L}(P)$,
\begin{align*}
A(x)&=A(\sqcup k(x))
=\bigvee_{y\in P}A(y)\wedge k(x)(y)\\
&=\bigvee_{y\in K(P)}A(y)\wedge {\uparrow} y(x).
\end{align*}
\noindent Thus, $\{{\uparrow} y\mid y\in K(P)\}$ is a base of $\Sigma_L P$. Since for every $y\in K(P)$,  ${\uparrow} y$ is a super-compact open set, $\Sigma_L P$ is  strong locally super-compact.
\hfill$\Box$

Conversely,  the following shows that every strong locally super-compact $L$-topological space can induce an algebraic $L$-dcpo.

\begin{tm}\label{tm-top-alg}  Let $X$ be a strong locally super-compact $L$-topological space. Then $({\rm pt}_{L}\mathcal{O}(X), {\rm sub}_{{O}(X)})$ is an algebraic $L$-dcpo.
\end{tm}
We first  prove two lemmas.
\begin{lm}\label{lm-bas-alg}Let $P$ be an  $L$-dcpo and $x\in P$. If there exists a directed L-subset $D$ such that $D\leq k(x)$ (cf. Definition \ref{dn-alg-dom} for the symbol $ k(x)$) and $\sqcup D=x$, then $ k(x)$ is directed and $x=\sqcup k(x)$.
\end{lm}
\noindent{\bf Proof.} It follows from $D\leq k(x)\leq {\downarrow} x$ that   $\sqcup k(x)=x$.
Next, we show  that $ k(x)$ is  directed.
Since $D$ is directed and $D\leq k(x)$,  we have $\bigvee_{y\in P} k(x)(y)=1$.
For every $y_{1}, y_2\in K(P)$,  we have
\begin{align*}
k(x)(y_1)\wedge k(x)(y_2)
&=e(y_1, x)\wedge e(y_2, x)\\
&=e(y_1, \sqcup D)\wedge e(y_2, \sqcup D)\\
&= (\bigvee_{d\in P}D(d)\wedge e(y_1, d))\wedge (\bigvee_{d\in P}D(d)\wedge e(y_2, d))\\
&=\bigvee_{d_1, d_2\in P}(D(d_1)\wedge D(d_2)\wedge e(y_1, d_1)\wedge e(y_2, d_2))\\
&=\bigvee_{d_1, d_2\in P}\bigvee_{d\in P}(D(d)\wedge e(d_1, d)\wedge e(d_2, d)\wedge e(y_1, d_1)\wedge e(y_2, d_2))\\
&\leq\bigvee_{d\in P} k(x)(d)\wedge e(y_1, d)\wedge e(y_2, d).
\end{align*}
Thus, $ k(x)$ is  directed.
\hfill$\Box$

\begin{lm}\label{lm2-th-alg} Let $X$ be an $L$-topological space and $p\in {\rm pt}_L\mathcal{O}(X)$. If $A$ is a super-compact open set, then ${\rm sub}_{\mathcal{O}(X)}([A], p)=p(A)$.
\end{lm}
\noindent{\bf Proof.} For every $B\in \mathcal{O}(X)$, by Proposition \ref{lm-pt}(1), we have ${\rm sub}_{X}(A, B)\wedge p(A)\leq  p(B)$.  Thus,
\vskip 3pt
\centerline{$p(A)\leq\bigwedge_{B\in \mathcal{O}(X)}{\rm sub}_{X}(A, B)\rightarrow p(B)={\rm sub}_{\mathcal{O}(X)}([A], p)$.}
\vskip 3pt
Since $A\in  \mathcal{O}(X)$, it follows that $\bigwedge_{B\in \mathcal{O}(X)}{\rm sub}_{X}(A, B)\rightarrow p(B)=p(A)$, as desired.
\hfill$\Box$

\noindent{\bf Proof of Theorem \ref{tm-top-alg}.} By Proposition \ref{lm-pt}(2), we only need to show that ${\rm pt}_L\mathcal{O}(X)$ is algebraic. Let $\mathcal{B}$ be a base consisting of super-compact open sets. For $p\in {\rm pt}_{L}\mathcal{O}(X)$, define $\mathbb{D}_p:{\rm pt}_{L}\mathcal{O}(X)\longrightarrow L$ as follows:
 $$
\ \ \mathbb{D}_p(q)=\left\{\begin{array}{ll}p(B),& q=[B], B\in \mathcal{B};\\
\  \ 0,& \mbox{otherwise}.
\end{array}\right.
$$
It is clear that $\mathbb{D}_p$ is well-defined.  We will prove that $\mathbb{D}_p$ is directed. Since $\mathcal{B}$ is a base, we have
$$\bigvee_{q\in {\rm pt}_{L}\mathcal{O}(X)}\mathbb{D}_p(q)=\bigvee_{B\in \mathcal{B}}p(B)=p(\bigvee_{B\in \mathcal{B}} B)=p(1_X)=1.$$
For every $p_1$, $p_2\in {\rm pt}_{L}\mathcal{O}(X)$,  without loss of generality, we assume that there exist $B_1, B_2\in \mathcal{B}$ such that $p_1=[B_1]$ and $p_2=[B_2]$. Then
\begin{align*}
\mathbb{D}_p(p_1)\wedge\mathbb{D}_p(p_2)&=p(B_1)\wedge p(B_2)
=p(B_1\wedge B_2)\\
&=p(\bigvee_{B\in \mathcal{B}}{\rm sub}(B, B_1\wedge B_2)\wedge B)\\
&=\bigvee_{B\in \mathcal{B}}{\rm sub}(B, B_1\wedge B_2)\wedge p(B)\\
&=\bigvee_{B\in \mathcal{B}}\mathbb{D}_p([B])\wedge {\rm sub}(B, B_1)\wedge {\rm sub}(B,  B_2)\\
&=\bigvee_{B\in \mathcal{B}}\mathbb{D}_p([B])\wedge {\rm sub}([B_1], [B])\wedge {\rm sub}([B_2], [B])\\
&=\bigvee_{q\in {\rm pt}_L\mathcal{O}(X)}\mathbb{D}_p(q)\wedge {\rm sub}(p_1, q)\wedge {\rm sub}(p_2, q).
\end{align*}
Therefore, $\mathbb{D}_p$ is  directed.
For every $A\in \mathcal{O}(X)$,
\begin{align*}
\sqcup \mathbb{D}_p(A)&=\bigvee_{q\in {\rm pt}_L\mathcal{O}(X)}\mathbb{D}_p(q)\wedge q(A)
=\bigvee_{B\in \mathcal{B}}\mathbb{D}_p([B])\wedge [B](A)\\
&=\bigvee_{B\in \mathcal{B}}p(B)\wedge {\rm sub}(B, A)
=p(\bigvee_{B\in \mathcal{B}}{\rm sub}(B, A)\wedge B)\\
&=p(A).
\end{align*}
By the arbitrariness of $A\in \mathcal{O}(X)$, we have $\sqcup \mathbb{D}_p=p$.

We will show that for every
$B\in \mathcal{B}$, $$\mathbb{D}_p([B])=k(p)([B]).$$
 For every $\mathbb{I}\in Idl_{L}({\rm pt}_L\mathcal{O}(X))$, by Lemma \ref{lm2-th-alg},
 \begin{align*}
 {\rm sub} ([B], \sqcup \mathbb{I})&={\rm sub}([B], \bigvee_{q\in {\rm pt}_L\mathcal{O}(X)} \mathbb{I}(q)\wedge q)
 =\bigvee_{q\in {\rm pt}_L\mathcal{O}(X)} \mathbb{I}(q)\wedge q(B)\\
 &=\bigvee_{q\in {\rm pt}_L\mathcal{O}(X)} \mathbb{I}(q)\wedge {\rm sub}([B], q)
 \leq \mathbb{I}([B]).
\end{align*}
By the arbitrariness of $\mathbb{I}\in Idl_{L}({\rm pt}_L\mathcal{O}(X))$, we have $[B]\in K({\rm pt}_L\mathcal{O}(X))$. Notice that
$$\mathbb{D}_p([B])=p(B)={\rm sub}([B], p).$$ We have $\mathbb{D}_p([B])= k(p)([B])$ (cf. Definition \ref{dn-alg-dom} for the symbol $ k(x)$). This shows that  $\mathbb{D}_p\leq k(p).$ It follows from Lemma \ref{lm-bas-alg} that $ k(p)$ is directed and $\sqcup k(p)=p$. By the arbitrariness of $p\in{\rm pt}_L\mathcal{O}(X)$, we have that ${\rm pt}_L\mathcal{O}(X)$ is an algebraic $L$-dcpo.\hfill$\Box$

\begin{rk} We now obtain  a representation for algebraic $L$-dcpo's. Specifically
 an $L$-ordered set $P$ is an algebraic $L$-dcpo if and only if there exists a strong locally super-compact  $L$-topological space $X$ such that $({\rm pt}_{L}\mathcal{O}(X), {\rm sub})\cong (P, e_p)$. Furthermore, if we impose the condition of  $L$-sobriety on  strong locally super-compact  $L$-topological spaces, we  obtain an isomorphism of  corresponding categories.
\end{rk}
Let $L\mathbf{\text{-}SLSCSob}$ denote the category of  strong locally super-compact sober spaces and continuous maps, and let $L\mathbf{\text{-}AlgDom}$ denote the category of  algebraic $L$-dcpo's and  Scott continuous maps. Clearly, $L\mathbf{\text{-}SLSCSob}$ and $L\mathbf{\text{-}AlgDom}$ are full subcategories of $L\mathbf{\text{-}LSCSob}$ and  $L\mathbf{\text{-}CDom}$, respectively. By Theorems \ref{tm-SO-id} and  \ref{tm-iso-sob-dom}, it follows that:

\begin{tm}\label{tm-iso-sob-alg} The categories $L\mathbf{\text{-}SLSCSob}$  and $L\mathbf{\text{-}AlgDom}$  are isomorphic.
\end{tm}

\section{Directed completions of continuous $L$-ordered sets}
As  an application of the results from the previous sections, in this section,  we study the  relationship between directed completions and $L$-sobrifications in the category of continuous  $L$-ordered sets.

Recall from Section 3 that the family
$\mathcal{O}{\rm pt}_{ L}\mathcal{O}(X)=\{\phi(A)\mid A\in \mathcal{O}(X)\}$
  is a $T_0$  $L$-topology, known as the spectral  $L$-topology on ${\rm pt}_{ L}\mathcal{O}(X)$. It is well-known that  ${\rm pt}_{ L}\mathcal{O}(X)$ equipped with the spectral $L$-topology is an $L$-sobrification of $X$. We often write ${\rm pt}_{ L}\mathcal{O}(X)$ with omitting the related  spectral $L$-topology.
 It should be noted that when we view ${\rm pt}_{ L}\mathcal{O}(X)$ as an $L$-ordered set, we always assume that it is equipped with the inclusion $L$-order. For simplicity, we often write ${\rm pt}_{ L}\mathcal{O}(X)$,  for $({\rm pt}_{ L}\mathcal{O}(X), {\rm sub}_{\mathcal{O}(X)})$.

 The following shows that the specialization $L$-order of $({\rm pt}_{ L}\mathcal{O}(X), \mathcal{O}{\rm pt}_{ L}\mathcal{O}(X))$ is precisely the  inclusion $L$-order.
\begin{pn}\label{pn-pt-spe} $\Omega_L {\rm pt}_{ L}\mathcal{O}(X)=({\rm pt_L}\mathcal{O}(X), {\rm sub}_{\mathcal{O}(X)})$.

\end{pn}
\noindent{\bf Proof.} For every pair $p$,$q$ in ${\rm pt}_{ L}\mathcal{O}(X)$, we have
\begin{align*}
e_{\mathcal{O}{\rm pt}_{ L}\mathcal{O}(X)}(p, q)&=\bigwedge_{A\in\mathcal{O}(X)}\phi(A)(p)\longrightarrow \phi(A)(q)\\
&=\bigwedge_{A\in\mathcal{O}(X)}p(A)\longrightarrow q(A)
={\rm sub}_{\mathcal{O}(X)}(p, q).
\end{align*}
Thus, $\Omega_L {\rm pt}_{ L}\mathcal{O}(X)=({\rm pt}_{ L}\mathcal{O}(X), {\rm sub}_{\mathcal{O}(X)})$.
\hfill$\Box$

By Proposition \ref{tm-dom-sob}, we know that the Scott $L$-topology of every  continuous $L$-ordered set is locally super-compact. The following shows that the locally super-compactness of a continuous $L$-ordered set is a topological property preserved under $L$-sobrification.

\begin{pn}\label{pn-fuzc-top}
Let $P$ be  a continuous $L$-ordered set. Then
 ${\rm pt}_{L}\sigma_{L}(P)$ is a locally super-compact $L$-sober space.
\end{pn}
\noindent{\bf Proof.} We only need to prove that ${\rm pt}_{L}\sigma_{L}(P)$  is locally super-compact. Let $A\in \sigma_{L}(P)$. For every $p\in {\rm pt}_L\sigma_{L}(P)$, by Proposition \ref{lm-uu-op}(3), we have
\begin{align*}
\phi(A)(p)
&=p(A)
=p(\bigvee_{a\in P}A(a)\wedge\uuar a)\\
&=\bigvee_{a\in P}A(a)\wedge \phi(\uuar a)(p)\\
&\leq \bigvee_{a\in P}{\rm sub}_{{\rm pt}_L\sigma_{L}(P)}({\uparrow} [a], \phi(A))\wedge({\uparrow} [a])^\circ(p).
\end{align*}
where ${\uparrow} [a](\mbox{-})={\rm sub}_{\sigma_{L}(P)}([a], \mbox{-})$.
The above  inequality relies on the following two  facts:
 \begin{itemize}
 \item$A(a)\leq {\rm sub}_{{\rm pt}_L\sigma_{L}(P)}({\uparrow} [a], \phi(A))$;

 \item $\phi (\uuar a)\leq ({\uparrow}[a])^{\circ}$.
 \end{itemize}

To prove the first fact, consider that  for every $q\in {\rm pt}_L\sigma_{L}(P)$, we have
$$A(a)\wedge {\rm sub}_{\sigma_{L}(P)}([a], q)\leq q(A)=\phi(A)(q).$$
Thus, $A(a)\leq \bigwedge_{q\in {\rm pt}_L\sigma_{L}(P)}{\rm sub}_{\sigma_{L}(P)}([a], q)\longrightarrow \phi(A)(q)$. This confirms the first fact.  Now, to show the second fact,
for every
$U\in \sigma_{L}(P)$ and $q\in {\rm pt}_L\sigma_{L}(P)$, by Proposition \ref{lm-pt}(1),
\vskip 6pt
\centerline{$q(\uuar a)\wedge U(a)\leq q(\uuar a)\wedge {\rm sub}_{P}(\uuar a, U)\leq q(U)$.}
By the arbitrariness of $U$, we have
\vskip 6pt
\centerline{$\phi(\uuar a)(q)\leq\bigwedge_{U\in \sigma_{L}(P)}[a](U)\rightarrow q(U)={\rm sub}_{\sigma_{L}(P)}([a], q)= {\uparrow}[a](q).$}
Thus,  $\phi(\uuar a)\leq {\uparrow}[a]$. Since $\phi(\uuar a)\in  \mathcal{O}{\rm pt}_L\sigma_{L}(P)$,  the second fact is verified.

 To sum up, we have
 $$\phi(A)=\bigvee_{a\in P}{\rm sub}_{{\rm pt}_L\sigma_{L}(P)}({\uparrow} [a], \phi(A))\wedge({\uparrow} [a])^\circ.$$
  Since every member of $\mathcal{O}{\rm pt}_L\sigma_{L}(P)$ is an upper set of $({\rm pt}_L\sigma_{L}(P), {\rm sub}_{\sigma_{L}(P)})$,  for every $a\in P$, ${\uparrow} [a]$ is a super-compact set. Therefore,  ${\rm pt}_{L}\sigma_{L}(P)$ is a locally super-compact space.\hfill$\Box$

\begin{pn}\label{pn-fuzc-top2}
Let $P$ be  a continuous $L$-ordered set. Then

 {\rm (1)} ${\rm pt}_{ L}\sigma_{L}(P)$ is a continuous $L$-dcpo.

{\rm (2)} $\sigma_{L}({\rm pt}_{L}\sigma_{L}(P))= \mathcal{O}{\rm pt}_{ L}\sigma_{L}(P).$ In other words,  the Scott $L$-topology on $({\rm pt}_{L}\sigma_{L}(P), {\rm sub}_{\sigma_{L}(P)})$ coincides with the the  related spectral $L$-topology.

\end{pn}
\noindent{\bf Proof.} (1) It  follows  directly from Propositions \ref{pn-sob-cont}, \ref{pn-pt-spe} and \ref{pn-fuzc-top}.

(2) By Proposition \ref{pn-pt-spe}, we have $\Omega_L {\rm pt}_{ L}\sigma_{L}(P)=({\rm pt}_{L}\sigma_{L}(P), {\rm sub})$. By Theorem \ref{tm-SO-id}(1), we have $\sigma_{L}({\rm pt}_{L}\sigma_{L}(P), {\rm sub})= \mathcal{O}({\rm pt}_{ L}(\sigma_{L}(P))).$
\hfill$\Box$

Inspired by the Xu and Mao's work (see \cite{Mao-Xu-2005}) and  Zhao's work (see \cite{zhao-fan}), we will introduce the notion of directed completions of continuous $L$-ordered sets.  However,  in \cite{Mao-Xu-2005}, the authors
 did not provide the notion of directed completion from a categorical perspective  to ensure the completion has the universal property in the sense of \cite{Bishop}.  In \cite{zhao-fan},  Zhao and Fan  studied a new type of  directed completion of posets,  called D-completion, and showed this  completion has the universal property. In \cite{Keimel-Lawson},  Keimel and Lawson   studied a general categorical construction via reflection functors for  the D-completion of $T_0$-spaces. Zhang et al. \cite{Zhang-Shi-Li} explored  the connection between the D-completion of posets and   D-completion (note that it is not a sobrification) of topologies in the lattice-valued setting, using the language of   closed sets within the framework of $L$-cotopology.  In the following, we formally provide a definition of directed completion applicable to continuous posets in a new way. Using the language of open sets, we will establish the connection between this type of directed completion and $L$-sobrification, ensuring that this connection can be extended to the frame-valued setting.

\begin{dn}
 Let $P$  be  a continuous $L$-ordered set and $Q$ be a continuous $L$-dcpo. If $j:P\longrightarrow Q$ is Scott continuous, then   $(Q, j)$, or $Q$, is called a directed completion of $P$ if for every continuous $L$-dcpo $M$ and Scott continuous map $f: P\longrightarrow M$, there exists a unique Scott continuous map $\overline{f}: Q\to M$ such that $\overline{f}\circ j=f$, i.e., the following diagram commutes.
\begin{displaymath}
\xymatrix@=8ex{P\ar[r]^{j}\ar[dr]_{f}&Q\ar@{-->}[d]^{\exists !\overline{f}}\\&M}
\end{displaymath}
\end{dn}
By the  universal property of directed completion,  directed completions of a continuous $L$-ordered set are unique up to $L$-order-isomorphism.

\begin{tm}\label{tm-dcpo-com}
Let $P$ be  a continuous $L$-ordered set. Then
$({\rm pt}_{L}\sigma_{L}(P), {\rm sub}_{\sigma_{L}(P)})$ with  $\eta_{P}$ is a directed completion of  $P$.
\end{tm}
\noindent{\bf Proof.}
By Proposition \ref{pn-fuzc-top2}(2), we have
 $\Sigma_L{\rm pt}_{L}\sigma_{L}(P)$ is an $L$-sobrification of $\Sigma_L P$.
By  Proposition \ref{pn-fuzc-top2}(1), we have $({\rm pt}_{L}\sigma_{L}(P), {\rm sub})$ is a continuous $L$-dcpo. Then it follows from  Proposition \ref{pn-scot-map}(2) that $\eta_{P}: P\longrightarrow {\rm pt}_{L}\sigma_{L}(P)$ is Scott continuous. For every continuous $L$-dcpo $Q$, and Scott continuous map $f:P\longrightarrow Q$, by Proposition \ref{pn-scot-map}(1) and Lemma \ref{lm-dom-sob}, we know that $\Sigma_L Q$ is $L$-sober and $f: \Sigma_L P\to \Sigma_L Q$ is continuous. Since $\Sigma_L{\rm pt}_{L}\sigma_{L}(P)$  is an $L$-sobrification of $\Sigma_L P$, we know that there exists a unique continuous map $\overline{f}: \Sigma_L{\rm pt}_{L}\sigma_{L}(P) \to \Sigma_L Q$ such that $\overline{f}\circ\eta_{P}=f$. By  Proposition \ref{pn-scot-map}(2),   $\overline{f}: {\rm pt}_{L}\sigma_{L}(P) \to Q$ is Scott continuous, as desired.\hfill$\Box$

Given a continuous poset $P$,  Xu and Mao \cite{Mao-Xu-2005} referred to the family of the irreducible Scott closed  sets (w.r.t  the set-theoretic inclusion order) as the directed completion of $P$. It is well-known that the set of all irreducible  Scott closed  sets is order-isomorphic to the set of  completely prime filters of the Scott open set lattice.  For a continuous
$L$-ordered set $P$, notice that every member of ${\rm pt}_{L}\sigma_{L}(P)$ is precisely the lattice-valued type of a completely prime filter.  By Theorem \ref{tm-dcpo-com},  the directed completion in the sense of Xu and Mao coincides with our's  directed completion  up to order-isomorphism. Therefore, Xu and Mao's directed completion also has the universal property.

Theorem \ref{tm-dcpo-com} also shows  that every continuous $L$-dcpo can be
retrieved from its Scott $L$-topology. Let $L\mathbf{\text{-}CPos}$ denote the category of all continuous $L$-ordered sets and Scott continuous maps. The above result shows that $L\mbox{-}{\bf CDom}$ is  reflective in $L\mbox{-}{\bf CPos}$. This result is a frame-valued type of Mislove's result (see \cite[Corollary 4.19]{Mislove}). Our approach differs from Mislove's. While Mislove utilizes both Scott open sets and Scott closed sets, we focus on the so-called points in terms of open sets.

In  the following, we provide  a characterization of  the directed completions of continuous $L$-ordered sets via $L$-sobrification.
\begin{tm}\label{tm-char-comp}
Let $P$ and $Q$ be two continuous $L$-ordered sets.

$(1)$  If $Y$ is an $L$-sobrification of $\Sigma_L P$, then  $\Omega_L Y$ is a directed completion of $P$;

$(2)$ If $Q$ is a directed completion of $P$, then  $\Sigma_L Q$ is an $L$-sobrification of $\Sigma_L P$;

$(3)$
  $Q$ is a directed completion  of $P$ if and only if  $\Sigma_L Q$ is an $L$-sobrification of $\Sigma_L P$.
\end{tm}
\noindent{\bf Proof.}
(1)  It is clear  that
$(Y, \mathcal{O}(Y))\cong ({\rm pt}_L\sigma_L (P), \mathcal{O}{\rm pt}_L\sigma_L (P)).$
Therefore, by Proposition \ref{pn-pt-spe}, we have
$\Omega_L Y\cong( {\rm pt}_L\sigma_L (P), {\rm sub}).$
It follows from Theorem \ref{tm-dcpo-com} that $\Omega_L Y$  is a directed completion of  $P$.

(2) It follows from Theorem \ref{tm-dcpo-com} that $(Q, e_Q) \cong( {\rm pt}_L\sigma_L (P), {\rm sub})$. Therefore $\Sigma_L Q\cong\Sigma_L{\rm pt}_L\sigma_L (P)$. It follows from Proposition \ref{pn-fuzc-top2}(2)
that $\Sigma_L Q$ is an $L$-sobrification of $\Sigma_L P$.

(3) By Proposition \ref{pn-SO-id},  we have $\Omega_L\Sigma_L (Q, e_Q)=(Q, e_Q)$.
Consequently,   Parts (1) and (2) implies that Part (3) holds.
\hfill$\Box$

\begin{rk}\label{rk-sob-D}

$(1)$ Given a continuous $L$-ordered set $P$, we now  have the following corresponding classes of directed completions  and  $L$-sobrifications:

\begin{picture}(0,35)
\put(189,18.5){$\Sigma_L$}\put(189,-2){$\Omega_L$}

\put(20,10){$\{ \mbox{directed completions of } P\}$}\put(175,15){\vector(1,0){37}}\put(212,10){\vector(-1,0){37}}
\put(220,10){$\{ L\mbox{-sobrifications of } (P, \sigma_{L}(P))\}$}
\end{picture}

Let ${\bf D}_L(P)$  denote the category of directed completions of  $P$ and let ${\bf S}_L(P)$ denote the category of $L$-sobrifications of  $(P, \sigma_{L}(P))$. These categories  are  full subcategories of $L\mathbf{\text{-}CDom}$ and
$L\mathbf{\text{-}LSCSob}$,   respectively.  By Theorems \ref{tm-SO-id} and \ref {tm-char-comp}, it is straightforward  to verify that the functors $\Sigma_L$ and $\Omega_L$ provide a categorical  isomorphism  between ${\bf D}_L(P)$  and ${\bf S}_L(P)$.

$(2)$ In the lattice-valued setting, for Theorem \ref{tm-char-comp}(3), the condition of  continuity of $Q$ is crucial  to ensure that $\Omega_L\Sigma_L (Q, e_Q)=(Q, e_Q)$ holds. However,  in    the classical setting,  this property  holds for every poset $Q$, regardless of whether
$Q$ is continuous. Here,  we indeed provide a topological characterization of directed completions of a continuous poset: a poset $Q$ is a directed completion  of a continuous poset $P$ if and only if  $(Q, \sigma(Q))$ is a sobrification of $(P, \sigma(P))$, where the symbol $\sigma(Q)$ denotes the Scott topology of the poset $Q$.

$(3)$ Based on  {\rm \cite[Remark 3]{zhao-fan}}, we deduce that  the D-completions   of a continuous poset $P$ in the sense of Zhao and   Fan {\em \cite{zhao-fan} } coincide with the   directed completions  of a continuous poset $P$ in our sense. By  {\em \cite[Theorem 5.5]{Zhang-Shi-Li}}, we deduce that a poset $Q$ is a D-completion  of a  poset $P$ if and only if  $(Q, \sigma(Q))$ is a D-completion of $(P, \sigma(P))$. Thus, for every continuous poset $P$, we summarize the following:
\begin{align*}
\ & Q \mbox{ is a directed completion of } P\\
\Longleftrightarrow\ &  Q \mbox{ is a D-completion of } P\\
\Longleftrightarrow\ & (Q, \sigma(Q)) \mbox{ is a D-completion of }(P, \sigma(P))\\
\Longleftrightarrow\ & (Q, \sigma(Q)) \mbox{ is a sobrification of }(P, \sigma(P)).
\end{align*}
\end{rk}
%
We end this section with two examples.
The first is a lattice-valued type of Lawson's example for round ideal completion (see \cite[Example 3.5]{Lawson}).
\begin{ex}\label{ex-roud}
Let $P$  be  a continuous $L$-ordered set. Define
  $${\bf RI}(P)=\{\bigvee_{x\in P}D(x)\wedge\dda x\mid D\in \mathcal{D}_{L}(P)\}.$$
Then $({\bf RI}(P), {\rm sub}_{P})\cong ({\rm pt}_{ L}\sigma_{L}(P), {\rm sub}_{    \sigma_L(P)})$. Thus,   $({\bf RI}(P), {\rm sub}_{   P})$ is a directed completion of  $P$ and $\Sigma_L({\bf RI}(P)
, {\rm sub}_P)$ is an $L$-sobrification of $\Sigma_L P$.
\end{ex}
\noindent{\bf Proof.} We divide the proof into three steps.

{\bf Step 1.} For every $p\in {\rm pt}_{ L}(\sigma_{L}(P))$,
put $D\in L^{P}$ by  $D(x)=p(\uuar x)$.
by Propoaitions \ref{lm-aux-way}, \ref{lm-con-int} and \ref{lm-uu-op},  it is routine  to verify that $D$ is directed.

Now  define a  map $f:{\rm pt}_{ L}\sigma_{L}(P)\to {\bf RI}(P)$ by
$$f(p)=\bigvee_{x\in P}p(\uuar x)\wedge \dda x.$$  It is clear  that $f$ is $L$-order-preserving.

 {\bf Step 2.} For $I\in {\bf RI}(P)$, there exists a directed $L$-subset $D_{I}$ such that $I=\bigvee_{x\in P}D_I(x)\wedge\dda x$. Define $\mathbb{D}\in L^{{\rm pt}_{ L} \sigma_{L}(P)}$ by
$$
\ \ \mathbb{D}(q)=\left\{\begin{array}{ll}I(x),& q=[x], x\in P;\\
\  \ 0,& \mbox{otherwise}.
\end{array}\right.
$$
Since $\Sigma_{L} P$ is $T_{0}$,  $\mathbb{D}$ is well-defined. We assert that $\mathbb{D}$ is directed. First, we have:
\begin{align*}
\bigvee_{x\in P}\mathbb{D}([x])&=\bigvee_{x\in P}I(x) =\bigvee_{x\in P}(\bigvee_{y\in P}D_{I}(y)\wedge \dda y)(x)\\
&=\bigvee_{y\in P}(D_{I}(y)\wedge \bigvee_{x\in P}\dda y(x))=\bigvee_{y\in P}D_{I}(y)=1.
\end{align*}
\noindent Second, for every $x_1, x_2\in P$, by the directedness of $\dda x$ $(\forall x\in P)$, we have
\begin{align*}
 \mathbb{D}([x_{1}])\wedge\mathbb{D}([x_{2}])&=I(x_1)\wedge I(x_2)\\
 &=\bigvee_{y_1, y_2\in P}D_{I}(y_1)\wedge D_{I}(y_2)\wedge \dda  y_1(x_1)\wedge \dda  y_2(x_2)\\
 &\leq\bigvee_{y_1, y_2, y_3\in P}D_{I}(y_3)\wedge e(y_1, y_3)\wedge e(y_2, y_3)\wedge \dda  y_1(x_1)\wedge\dda  y_2(x_2)\\
&\leq\bigvee_{ y_3\in P}D_{I}(y_3)\wedge \dda y_3(x_1)\wedge \dda  y_3(x_2)\\
&\leq\bigvee_{ y_3\in P}D_{I}(y_3)\wedge \bigvee_{x_3\in X}(\dda y_3(x_3)\wedge e(x_1, x_3)\wedge e(x_2, x_3))\\
&=\bigvee_{x_3\in P}(\bigvee_{ y_3\in P}D_{I}(y_3)\wedge\dda y_3(x_3))\wedge e(x_1, x_3)\wedge e(x_2, x_3)\\
&=\bigvee_{ x_3\in P} I(x_3)\wedge e(x_1, x_3)\wedge e(x_2, x_3)\\
&\leq \bigvee_{ x_3\in P} \mathbb{D}([x_{3}])\wedge {\rm sub}([x_1], [x_3])\wedge {\rm sub}([x_2], [x_3]).
\end{align*}
 Thus, $\mathbb{D}$ is directed. By Proposition \ref{lm-pt}(2),  we have $\sqcup \mathbb{D}=\bigvee_{x\in P}I(x)\wedge [x]$. Hence,  $\bigvee_{x\in P}I(x)\wedge [x]\in {\rm pt}_{ L}\sigma_{L}(P)$. Now,  define a map $g:{\bf RI}(P)\to {\rm pt}_{ L}\sigma_{L}(P)$ by
 $$g(I)=\bigvee_{x\in P}I(x)\wedge [x].$$ It is clear that $g$ is $L$-order-preserving.

 {\bf Step 3.}  We need to show that $f$ and $g$ are inverse to each other. For every $p\in {\rm pt}_{ L}\sigma_{L}(P)$ and $A\in \sigma_{L}(P)$, by Proposition \ref{lm-uu-op}(3), we have:
\begin{align*}
p(A)&=p(\bigvee_{x\in P}A(x)\wedge \uuar x)=\bigvee_{x\in P}A(x)\wedge p(\uuar x)
=\bigvee_{x\in P}(\bigvee_{y\in P}A(y)\wedge\uuar y(x))\wedge p(\uuar x)\\
&=\bigvee_{y\in P}(\bigvee_{x\in P} p(\uuar x)\wedge \dda x(y))\wedge [y](A)
=\bigvee_{y\in P}(f(p)(y))\wedge [y](A)
=g(f(p))(A).
\end{align*}
Thus, $gf(p)=p$.

For every $I=\bigvee_{y\in P}D(y)\wedge \dda y\in {\bf RI}(P)$, for every $t\in P$, by Proposition \ref{lm-con-int},
\begin{align*}
\bigvee_{x\in P}I(x)\wedge \dda x(t)&=\bigvee_{x\in P}(\bigvee_{y\in P}D(y)\wedge \dda y(x))\wedge \dda x(t)
=\bigvee_{y\in P}D(y)\wedge \dda y(t)=I(t).
\end{align*}
Therefore,  $I=\bigvee_{x\in P}I(x)\wedge \dda x$.
For every $t\in P$,  we have:
\begin{align*}
f(g(I))(t)&=f(\bigvee_{x\in P}I(x)\wedge [x])(t)
=\bigvee_{x, y\in P}I(x)\wedge \uuar y(x)\wedge \dda y(t)\\
&=\bigvee_{x\in P}I(x)\wedge \dda x(t)=I(t).
\end{align*}
Thus, $fg(I)=I$. \hfill$\Box$

The above Example shows that, for a continuous poset, its  round ideal completion coincides with its directed completion  up to order-isomorphism.

\begin{ex}\label{ex-rou-ide}
Let $X=\{x, y\}$ and $L=\{0,a,b,1\}$ be a frame with $0\leq a, b\leq1$ and $a\parallel b$. Define an $L$-order $e$ on $X$ by setting
$e(x,x)=1$, $e(x, y)=1$, $e(y, y)=1$ and $e(y, x)=0$. The ideals of $(X, e)$ and their relative supremums, if they exist, are shown in the following table.

$$\begin{tabular}{|c|c|c|c|c|}
\hline
  &$\ I_1\ $&$\ I_2\ $&$\ I_3\ $&$\ I_4\ $                      \\
\hline
$x$ &$ 1$&$1$&$ 1$&$1$   \\
\hline
$y$ &$ 1$&$a$&$ b$&$0$   \\
\hline
{\rm sup} &$ y$&$\times$&$\times$&$x$   \\
\hline
\end{tabular}
$$

We observe that  $(X, e)$ is not an $L$-dcpo. By  a routine calculation, we have the following table.
$$\begin{tabular}{|c|c|c|c|}
\hline
  &$\ x\ $&$\ y\ $&{\rm sup}                     \\
\hline
$\dda y$ &$ 1$&$1$&$ y$   \\
\hline
$\dda x$ &$ 1$&$0$&$ x$  \\
\hline

\end{tabular}
$$

 Thus, $(X, e)$ is a continuous $L$-ordered set, but not a continuous $L$-dcpo.
 The Scott open sets of $(X, e)$ are illustrated in Figure 6.1,  where each pair $(m, n)$ denote the $L$-subset with  value $m$ on $y$ and  value $n$ on $x$.

\begin{center}
  \includegraphics[height=1.6in,width=3.557in]{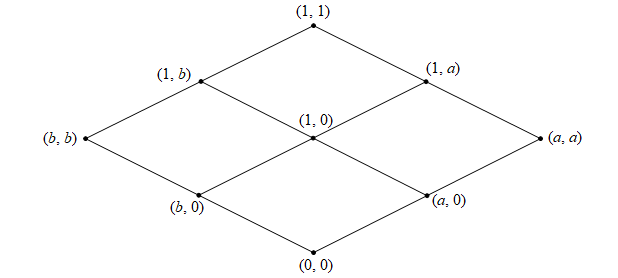}\\
 {\scriptsize {\rm {\bf Fig. 6.1:} The Scott $L$-topology on $(X,e)$} }
 \end{center}

Through computation, the Scott open sets of  $(L, e_L)$ are depicted in  Figure 6.2,  where each tuple $(m, n, p, q)$ represents  the $L$-subset with  values $m, n, p$, and $q$ assigned to  $1, a, b$, and $0$, respectively.

\begin{center}
  \includegraphics[height=1.6in,width=3.557in]{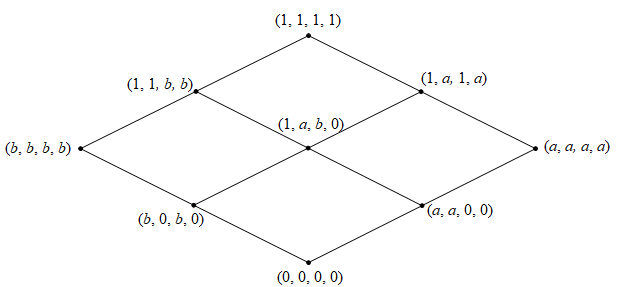}\\
 {\scriptsize {\rm {\bf Fig. 6.2:} The Scott $L$-topology on $(L, e_L)$} }
 \end{center}

By combining Figures 6.1 and 6.2 and performing mechanical calculations, we obtain  $(\sigma_{L}(L), {\rm sub}_{L})\cong (\sigma_{L}(X), {\rm sub}_{X})$. According to  Theorem \ref{tm-char-sobri}(3),  $(L, \sigma_{L}(L))$ is an $L$-sobrification of $(X, \sigma_{L}(X))$. Consequently, by Theorem \ref{tm-char-comp}, $(L, e_L)$ is a directed completion of $(X, e)$.

\end{ex}

\section{Conclusion}

This paper studies the topological representations of frame-valued domains via   $L$-sobriety for a frame $L$ as the truth value table. We prove that   the category of continuous $L$-dcpo's with Scott continuous maps is isomorphic to that of locally super-compact $L$-sober spaces with  continuous maps;
the category of algebraic $L$-dcpo's with Scott continuous maps is isomorphic to that of strong locally super-compact $L$-sober spaces with continuous maps. One byproduct of our  topological representation of frame-valued domains is the topological characterization of directed completions of  continuous $L$-ordered sets. This work also demonstrates that, just as sobriety plays an important role in domain theory/order theory, frame-valued sobriety introduced by \cite{De-Zhang} is essential in the study of frame-valued order theory. In summary, this research opens a way for finding  topological representations of frame-valued domains.


Our work was mainly inspired by Xu and Mao's work \cite{Mao-Xu-2005,XU-mao-form},  Mislove's work  \cite{Mislove}) and Lawson's work \cite{Lawson}. However, we would like to stress the fact that the results of this paper are not merely an exercise in fuzzyfication, as we did not follow their proof methods step by step.
 In order to generalize their findings to the  the frame-valued setting, a more direct proof was necessary. Additionally, we explored some results are  even new in the  classical setting and constructed some non-ordinary  examples to illustrate our findings.  The primary innovations of our work are as follows:

\begin{itemize}
\item In \cite{XU-mao-form}, Xu and Mao provided a representation for various  domains via the so-called formal points related to  super-compact quasi-bases of  sober spaces. We have replaced  the  formal points with the points (equivalently, completely prime filters) of the  open set  lattice.  These  replacements make the corresponding proofs more direct and allow us  to smoothly generalized the corresponding results to those in  the frame-valued setting.

\item In \cite{XU-mao-form}, Xu and Mao  did not explore the relationship between  domains  and locally super-compact sober spaces from the categorical perspective, while we  established the isomorphisms of relevant categories in this paper (see Theorems \ref{tm-iso-sob-dom}, \ref{tm-iso-sob-alg}).  In  obtaining categorical results, Proposition \ref{pn-scot-map} plays a key role,  which is well-known in the classical case. But in the frame-valued case,  its proof is new which relies on  the continuity of  the way-below relation.

\item In \cite{Mislove},  Mislove jointly  used open sets and closed sets  to show that the category of continuous dcpos is  reflective in that   of continuous posets. However,  in the frame-valued context, the table of truth-values does not satisfy the law of double negation, and there is no natural way to switch between open sets and closed sets.   So different ideas and techniques are needed to make the new proof feasible in the frame-valued setting.  Here, by means of the locally super-compactness of spectral topologies (see Proposition \ref{pn-fuzc-top}), we  deduced the corresponding result (see Theorem \ref{tm-dcpo-com}) in the language of open sets. This proof method is new and  can be considered as an application of Theorem \ref{tm-iso-sob-dom}.

\item The notion of directed completion of continuous posets appeared in Xu and Mao's work \cite{Mao-Xu-2005}. Since it is customary to have the result that the corresponding completion of $P$ should have the universal property in the sense of \cite{Bishop}, we present this notion  in a new way to ensure that the resulting directed completion has the universal property. By means of  the locally super-compacteness of spectral spaces, we derived Theorem \ref{tm-char-comp} which sheds light on the relationship between directed completion and sobrification.  Specifically, this is also a new result in the classical setting (see Remark \ref{rk-sob-D}). Moreover,  Example \ref{ex-roud} shows that, for a continuous poset, its  round ideal completion coincides with its directed completion  \cite{R. Hofmann}. In fact, these connections were not fully  revealed in \cite{R. Hofmann,Mislove,Mao-Xu-2005}.
\end{itemize}

\section*{Acknowledgment}

This paper is supported by National Natural Science Foundation of China
(12231007, 12371462), Jiangsu Provincial Innovative and Entrepreneurial
Talent Support Plan (JSSCRC2021521).

\end{document}